\documentclass[reprint, amsmath, amssymb, aps, pre, floatfix, showkeys, superscriptaddress]{revtex4-2}

\usepackage{mathrsfs}
\usepackage{graphicx}
\usepackage{dcolumn} 
\usepackage{bm}
\usepackage{mathtools}
\usepackage{color}
\usepackage{braket}
\usepackage{tikz}
\usepackage{pgfplots}
\usepackage{comment}
\pgfplotsset{compat=1.18}

\begin{document}

\title{Ghost states underlying spatial and temporal patterns: how non-existing invariant solutions control nonlinear dynamics}

\author{Zheng Zheng}
\author{Pierre Beck}
\author{Tian Yang}
\author{Omid Ashtari}
\affiliation{Emergent Complexity in Physical Systems Laboratory (ECPS), École Polytechnique Fédérale de Lausanne, CH-1015 Lausanne, Switzerland}
\author{Jeremy P Parker}\email{jparker002@dundee.ac.uk}
\affiliation{Emergent Complexity in Physical Systems Laboratory (ECPS), École Polytechnique Fédérale de Lausanne, CH-1015 Lausanne, Switzerland}
\affiliation{Department of Mathematics, University of Dundee, Dundee DD1 4HN, United Kingdom}
\author{Tobias M Schneider}\email{tobias.schneider@epfl.ch}
\affiliation{Emergent Complexity in Physical Systems Laboratory (ECPS), École Polytechnique Fédérale de Lausanne, CH-1015 Lausanne, Switzerland}
\date{\today}

\begin{abstract}
Close to a saddle-node bifurcation, when two invariant solutions collide and disappear, the behavior of a dynamical system can closely resemble that of a solution which is no longer present at the chosen parameter value. For bifurcating equilibria in low-dimensional ODEs, the influence of such `ghosts' on the temporal behavior of the system, namely delayed transitions, has been studied previously. 
We consider spatio-temporal PDEs and characterize the phenomenon of ghosts by defining representative state-space structures, which we term `ghost states,' as minima of appropriately chosen cost functions. Using recently developed variational methods, we can compute and parametrically continue ghost states of equilibria, periodic orbits, and other invariant solutions. We demonstrate the relevance of ghost states to the observed dynamics in various nonlinear systems including chaotic maps, the Lorenz ODE system, the spatio-temporally chaotic Kuramoto--Sivashinsky PDE, the buckling of an elastic arc, and 3D Rayleigh--B\'enard convection.
\end{abstract}

\maketitle

\section{Introduction}
\label{sec:intro}
Characterizing dynamical systems from a geometrical point of view provides a powerful approach to complex nonlinear systems, where obtaining analytical or closed-form solutions is typically impossible. In this approach, the state of the system at a given time is viewed as a `point,' and its evolution is viewed as a `trajectory' within an abstract, possibly very high-dimensional, and often formally infinite-dimensional state space. Instead of describing the behavior of individual solutions for given initial conditions, a geometrical approach characterizes the dynamics in terms of the topology of the ensemble of all these trajectories, known as the phase portrait of the system. The idea for such a geometrical approach was initiated in the late nineteenth century \citep{Poincare1892, Birkhoff1927}, which later became a key element in the study of deterministic chaos \citep{lorenz1963deterministic}. Thanks to advances in computational resources and methods over the past three decades, this approach has been successfully transferred from low-dimensional dynamical systems governed by ordinary differential equations (ODEs) to those governed by partial differential equations (PDEs). The function space of a PDE is infinite-dimensional, which, in discretized form, results in finite but very high-dimensional problems. Examples of such PDE systems include fluid flows \citep{Gibson2008, Kawahara2012} and nonlinear optics \citep{Tlidi1997, Parra-Rivas2021}, among others, where a geometrical approach has proven to be particularly useful.

A geometrical approach aims first to characterize how general trajectories are organized within the state space of a dynamical system, and then to discover how this organization changes as the control parameters of the system are varied. The shape of the general trajectories is determined by special trajectories corresponding to dynamically invariant solutions with a simple dependence on time. The simplest of these invariant solutions include steady states, where the system does not evolve over time, and periodic orbits, where the state recurs exactly after a given period. These invariant solutions might be dynamically unstable and, thus, never realized naturally in a numerical simulation or a laboratory experiment. However, their local dynamics and the dynamical connections between them dictate the shape of the other trajectories.

When the control parameters of a dynamical system are varied, a geometrical description is typically concerned with bifurcations. A bifurcation is the appearance of a topologically nonequivalent phase portrait at precise parameter values \citep{Kuznetsov2004}. This topological change involves the appearance or disappearance of fixed points, periodic orbits, or more complex invariant sets, or changes in their stability properties. A bifurcation is local if its associated topological change can be detected within any small neighborhood of the invariant solution involved in the bifurcation; if this change cannot be identified within small neighborhoods of the invariant solution, then the bifurcation is global. Through a local bifurcation, the shape and dynamical properties of trajectories outside a small neighborhood around the involved invariant solution do not change immediately. Consequently, whenever two invariant solutions collide and annihilate in a saddle-node (fold) bifurcation, the properties of these solutions are still `felt' by the dynamics even though the solutions are not present anymore. This phenomenon is called the `ghost' of the saddle-node bifurcation \citep{strogatz2018nonlinear}.

When two fixed points collide in a saddle-node bifurcation, slow evolution emerges in the ghost region as the key remaining property of the disappeared equilibria. The ghost region attracts and slows down the flow of nearby trajectories, delaying their transition to another attracting region in the state space. As the control parameter approaches the bifurcation value, the time spent in the vicinity of the ghost (or, equivalently, the transition time needed for the system to pass the ghost and move to some other equilibrium state) increases and diverges in a predictable manner. The scaling of the transition time with the distance from the bifurcation point in the parameter space has been extensively studied before \citep{Strogatz1989, Fontich2008, sardanyes2020noise}. Some authors \citep{Fontich2008, Canela2021}, particularly in quantum mechanics \citep{kus1993prebifurcation, sundaram1995standard}, interpret ghosts by extending the real-valued state space of the dynamics to the space of complex-valued state variables. In this picture, two fixed points in the real-valued subspace collide and disappear from that subspace, giving rise to two repelling fixed points outside the real-valued subspace. The effect of these complex-valued repellers on the real-valued subspace, i.e., the original state space of the system, is what we identify as the ghost phenomenon. Such an approach, however, does not generalize to very complicated systems such as PDEs.

Ghosts of equilibrium solutions have been investigated across a variety of applications. These include electronic circuits \citep{Trickey1998}, population dynamics and catalytic hypercycles of biochemical reactions \citep{sardanyes2006ghosts}, the buckling of elastic solid structures \citep{gomez2017critical}, and understanding and predicting financial crises \citep{smug2018predicting, diks2019critical}, to name a few. Ghosts of equilibria play a significant role in so-called critical transitions: when the control parameter of a system, which is at a stable equilibrium state, is changed slowly toward and then beyond a critical value, where the equilibrium solution disappears through a saddle-node bifurcation or loses stability through other bifurcations, the system transitions to a new attracting region of the state space. In the case of a saddle-node bifurcation, the transition time is controlled by the ghost. This phenomenon is known as a critical transition, and the critical threshold is often referred to as the tipping point. The significance of a critical transition is that the new equilibrium state may be undesirable while, due to strong hysteresis, the transition cannot be reversed immediately by returning to a parameter value for which the desirable equilibrium exists again \citep{scheffer2009early, Sardanyes2024}. Recently, \citet{calsina2023ghost} considered a one-dimensional (1D) PDE model for a spatially extended population dynamics. They studied how the strength of the diffusion term in their model affects the delayed transition times to extinction. The transition time is governed by the ghost, as the parameter value is chosen close to a saddle-node bifurcation through which two other equilibria disappear and leave the uniform zero as the only attracting solution. However, the majority of previous research has investigated simple ODE systems.

In this work, we formalize the ghost phenomenon by defining representative state-space structures without the necessity of being asymptotically close to the saddle-node bifurcation point in the parameter space. We term these representative structures the `ghost states.' Characterizing the ghost phenomenon in terms of the representative ghost states enables us to follow three main objectives that have not been fully addressed previously. First, unlike previous works that study the ghost phenomenon as the system approaches the bifurcation point in the parameter space, we are interested in this phenomenon as the system moves away from it, so that less time is expected to be spent in the ghost region. The ghost states enable us to investigate how the ghost region evolves and whether it remains relevant to the dynamics as the control parameter is varied further from the bifurcation value. Secondly, we are interested in the ghost phenomenon in spatially extended dynamical systems governed by PDEs, such as fluid flows. We aim to explore the relevance of ghosts to both the spatial and temporal properties of the dynamics, rather than solely the temporal characteristics of delayed transitions, which have been the primary focus of previous studies. If the numerical computation of the ghost states can be scaled to high-dimensional problems, as will be demonstrated, the proposed characterization of the ghost phenomenon can be applied to high-dimensional discretizations of PDEs as well as low-dimensional ODEs. Finally, in addition to equilibrium solutions, we are interested in the ghosts resulting from the saddle-node bifurcation of time-varying invariant solutions, which has not been explored previously. We propose a family of methods that, based on a single unifying idea, defines and computes the ghost states for different types of invariant solutions. This enables us to pursue the previous two objectives for the ghosts of invariant solutions of various types. Specifically, we study the ghosts of periodic orbits as well as equilibria in the present work.

We define the ghost state in the space of all sets that have the same topological structure as the invariant solutions involved in the saddle-node bifurcation. For the ghosts of equilibria and periodic orbits, this space includes all points and all loops in the state space, respectively. Within the respective search space, we formulate invariant solutions as the global minima (zeros) of a suitably defined non-negative cost function. Following this formulation, we define the ghost states as the non-zero minima of the cost function, which lift from zero as soon as the control parameter passes the bifurcation value. As the control parameter approaches the bifurcation value, the minimum value of the cost function decreases to zero, and the ghost state becomes the invariant solution itself, by construction and as expected. We demonstrate that as the control parameter is varied further from the bifurcation value, the ghost state captures the essential properties of the ghost phenomenon. Defined as such, the ghost of an equilibrium solution represents the locally slowest point in the state space at the chosen parameter value, and the ghost of a periodic orbit represents the best-fit loop to the vector field induced by the governing equations.
A recently developed family of variational methods \cite{farazmand2016adjoint,azimi2022constructing,parker2022variational,ashtari2023identifying} enables us to apply this approach to high-dimensional dynamical systems, including discretizations of PDEs.
However, the cost functions whose minima represent the ghost states are not uniquely defined. Therefore, the ghost states are not purely a property of the dynamical system, and their exact properties depend on the specific choice of the cost function. Nevertheless, we will see that the ghost states defined based on different cost functions consistently capture the characteristics of the ghost phenomenon, providing important insight into the dynamics of a system.

One particular application of ghosts is to intermittent chaos. In some instances (the type-1 transition of \citet{pomeau1980intermittent}), intermittent chaos is preceded by the saddle-node bifurcation of a periodic orbit. The trajectories of the subsequent chaotic attractor intermittently visit the region of the state space where the solutions formerly existed---the ghost of the bifurcation. We compute the ghost of such a periodic orbit in the Lorenz system. We demonstrate that, compared to the periodic orbit at the bifurcation point, the ghost state provides a better match with the regions of the state space frequently visited by the intermittent chaos at the studied parameter value.

This paper proceeds as follows: in Sec. \ref{sec:background}, we present the mathematical background of ghosts. Section \ref{sec:method} details the general principles of the family of methods we have developed. In Sec. \ref{sec:examples}, we present several examples of applying our methods to problems of various complexities, to show the power but also the pitfalls of our methods. In Sec. \ref{sec:ilc}, we apply our methods to the three-dimensional (3D) Rayleigh--B\'enard convection problem. We give concluding remarks in Sec. \ref{sec:conclusion}. 

\section{Background}
\label{sec:background}
Bifurcations occur in dynamical systems when invariant solutions, such as steady states and periodic orbits, appear, cease to exist, or change in stability properties at precise parameter values. In most of the usual classes of bifurcations, such as pitchfork, Hopf and transcritical bifurcations, invariant solutions exist on both sides of the bifurcation point in the parameter space, governing the dynamics of the system locally in the state space. By contrast, in saddle-node (fold) bifurcations, two invariant solutions collide and annihilate one another. Therefore, on one side of the bifurcation point in the parameter space, there is no longer an invariant solution in the state-space locality where the two solutions collide. However, close to this bifurcation point in the parameter space, the dynamics still `feel' the properties of the invariant solutions that are no longer present. This phenomenon is called the `ghost' of the saddle-node bifurcation. In this paper, we interchangeably use `the ghost of a bifurcation' or `the ghost of a solution.' Note that two invariant solutions, one of which has an additional unstable direction, collide at the bifurcation point. Therefore, which of the two solutions we are referring to is immaterial.

For lack of better terms, we shall assume in this paper, without loss of generality, that as the bifurcation parameter is \textit{increased}, a pair of invariant solutions collide in a saddle-node bifurcation. Therefore, we may refer to the region of the parameter space with the two solutions as being \textit{below} the critical value of bifurcation, and the region without the solutions as being \textit{above} it.

\begin{figure}
    \centering
    \includegraphics[width = \linewidth]{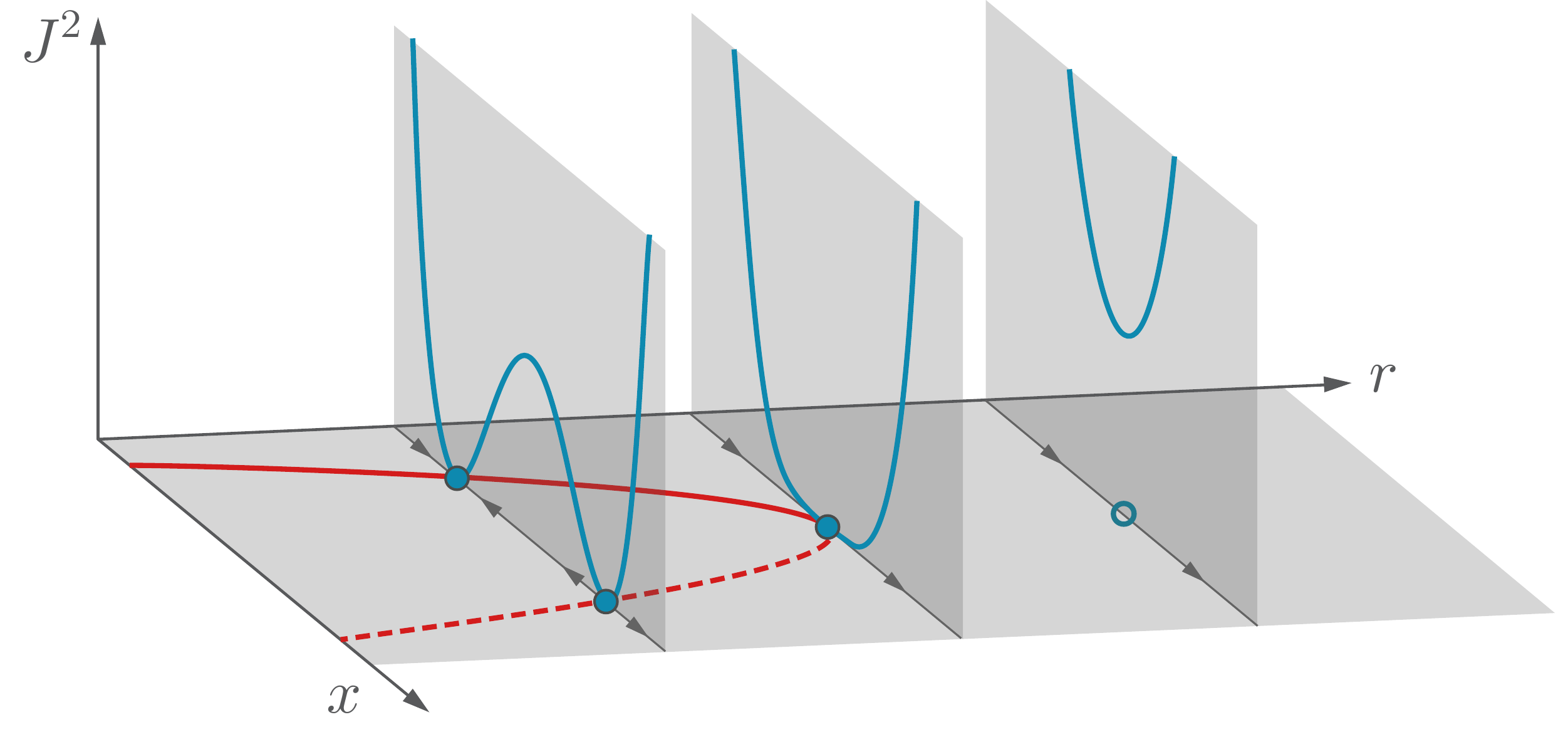}
    \caption{The 1D dynamical system given by Eq.~\eqref{eq:1DDS}. The red curve shows the locus of equilibrium states in the $r$--$x$ plane, undergoing a saddle-node bifurcation. For $r<0$ (left section), two equilibrium points exist, one stable and one unstable. When $r=0$ (middle), one equilibrium exists at $x=0$ which is neither stable nor unstable. For $r>0$ (right), no equilibria exist and the state moves in the positive $x$-direction for all times. A cost function $J^2 := (1/2)(r+x^2)^2$ penalizes the deviation from equilibrium such that $J=0$ at equilibrium points (filled markers) and $J>0$ otherwise. For the control parameter above the critical value of bifurcation, $r>0$, the ghost is defined as the state at which the cost function takes a non-zero minimum value (the open marker).}
    \label{fig:1DDS}
\end{figure}

\begin{figure}
    \centering
    \includegraphics[width = \columnwidth]{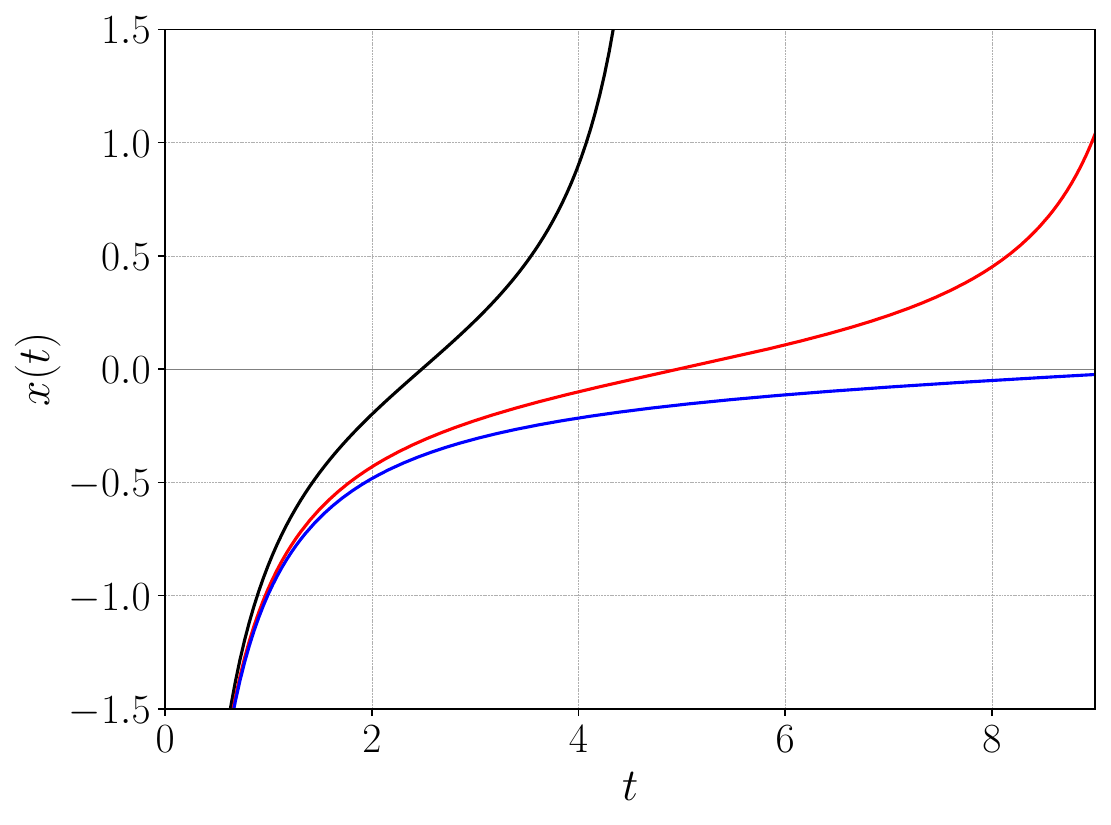}
    \caption{Solutions of $\mathrm{d}x/\mathrm{d}t=r+x^2$ with $x(0)=-\infty$ for different values of $r>0$. Black, red, and blue curves correspond to $r=0.4$, $r=0.1$, and $r=0.025$, respectively. Over time, $x$ monotonically increases toward infinity.
    However, there is an increasingly pronounced slowdown around $x=0$ as $r\to 0$.
    This delayed transition is the key temporal property of the ghost of equilibria past a saddle-node bifurcation.}
    \label{fig:tan}
\end{figure}

The canonical saddle-node bifurcation is found in the 1D real-valued ODE system \citep{strogatz2018nonlinear}
\begin{equation}
    \dfrac{\mathrm{d}x}{\mathrm{d}t}=r+x^2.
    \label{eq:1DDS}
\end{equation}
Immediately, we see that when $r<0$, there are two points which satisfy ${\mathrm{d}x}/{\mathrm{d}t}=0$, at $x=\pm\sqrt{-r}$. When $r=0$, we have a degenerate case with only one equilibrium at the origin. For $r>0$, no equilibrium exists at all; hence, $x$ monotonically increases toward infinity over time for all initial conditions. Fixed points and trajectories of this system are shown in Fig.~\ref{fig:1DDS}. However, this is not the full picture, as becomes apparent when we explicitly solve the equation for the $r>0$ case:
\begin{equation}
    x(t) = \sqrt{r} \tan{\left(t\sqrt{r}-\frac{\pi}{2}\right)},
\end{equation}
for $t\in (0,\pi/\sqrt{r})$,
where this solution satisfies $x\to-\infty$ as $t\to 0$ and $x\to\infty$ as $t \to \pi/\sqrt{r}$. This transition time diverges as $r\to 0$, with increasingly large amounts of time spent close to $x=0$, as shown in Fig.~\ref{fig:tan}. The slow evolution around $x=0$ is the remaining property of the non-existing equilibria in the ghost region. This effect is called critical slowdown, characterized by the well-known $r^{-1/2}$ scaling of the transition time \cite{Strogatz1989}. 

This simple example shows the phenomenon of ghosts for saddle-node bifurcations of steady states in continuous-time dynamical systems. However, the same concept extends to fixed points in maps, as will be demonstrated in Sec. \ref{sec:map}. Furthermore, we should expect analogous behavior in higher-dimensional invariant solutions such as periodic orbits and invariant tori, since these can be reduced to fixed points of maps via the use of Poincar\'e sections.

In the following, we characterize equilibria and periodic orbits as the absolute minima (zeros) of a suitably defined non-negative cost function. This enables us to define the ghost as a state-space structure for which the respective cost function takes a minimum yet non-zero value when the control parameter passes the saddle-node bifurcation value and lies above it.
This is shown for the dynamics \eqref{eq:1DDS} in Fig.~\ref{fig:1DDS}.

\section{Method}
\label{sec:method}
Equilibrium, periodic and quasi-periodic solutions are time-invariant sets of certain topological properties in the state space of the system. Equilibrium solutions are represented by isolated fixed points, periodic orbits by closed loops, and quasi-periodic orbits by higher-dimensional torus-shaped manifolds in the state space. We view the computation of an invariant solution as a minimization problem in the space of all sets possessing the same topological structure as the sought-after solution. In this approach, the deviation of such a set from satisfying the definition of the objective invariant solution is penalized by a non-negative cost function. The minimization of the cost function evolves a guess until at a global minimum, where the cost function takes the value of zero, an invariant solution is found. Such a cost function is not unique. In the following, a simple choice for equilibria and periodic orbits is presented in a general setting.

Suppose the following autonomous continuous-time dynamical system is given:
\begin{equation}
\label{eq:general_governing_equation}
    \dfrac{\partial\mathbf{u}}{\partial t}=\mathbf{f}(\mathbf{u}),
\end{equation}
where $\mathbf{u}(\mathbf{x},t)\in\mathbb{R}^n$ is a real-valued field defined over the spatial domain $\mathbf{x}\in\Omega\subseteq\mathbb{R}^d$ and time $t\in\mathbb{R}$, and $\mathbf{f}$ is a smooth nonlinear operator enforcing the boundary conditions (BCs) at $\partial\Omega$, the boundaries of $\Omega$.
Equilibria, that are time-independent solutions satisfying $\mathbf{f}(\mathbf{u})=\mathbf{0}$, can be found by minimizing
\begin{equation}
    \label{eq:ghosts_cost_FP}
    J_e^2 := \dfrac{1}{2}\int_\Omega\left\|\mathbf{f}(\mathbf{u})\right\|^2\,\mathrm{d}\mathbf{x},
\end{equation}
over all feasible states in the state space
\begin{equation}
    \label{eq:search_space_EQ}
    \mathscr{C}_e := \left\{\mathbf{u}(\mathbf{x})
    \left|
    \begin{array}{l}
        \mathbf{u}:\Omega\to\mathbb{R}^n\;\mathrm{sufficiently\; smooth}\\
        \mathbf{u}\mathrm{\;satisfies\;BCs\;at\;}\partial\Omega
    \end{array}
    \right.
    \right\}.
\end{equation}
Here, $\|\cdot\|$ denotes the standard Euclidean $L_2$-norm in $\mathbb{R}^n$. The cost function $J_e^2$ is zero if and only if $\mathbf{u}\in\mathscr{C}_e$ is an equilibrium solution of Eq.~\eqref{eq:general_governing_equation}, and takes a positive value otherwise. In a discrete-time system, an equivalent cost function is proportional to the squared $L_2$-distance between a state and its first iterate (see Sec. \ref{sec:map}).

In the vector field induced by the governing equation \eqref{eq:general_governing_equation}, periodic orbits are closed integral curves traversed in a finite time period. Hence, periodic orbits can be found by minimizing
\begin{equation}
    \label{eq:ghosts_cost_PO}
    J_p^2 := \dfrac{1}{2}\int_0^1\int_\Omega\left\|\dfrac{1}{T}\dfrac{\partial\mathbf{u}}{\partial s}-\mathbf{f}\left(\mathbf{u}\right)\right\|^2\,\mathrm{d}\mathbf{x}\mathrm{d}s,
\end{equation}
over the space of loops in state space
\begin{equation}
\label{eq:search_space_PO}
    \mathscr{C}_p := \left\{
    \begin{bmatrix}
        \mathbf{u}(\mathbf{x},s)\\
        T
    \end{bmatrix}\left|
    \begin{array}{l}
        \mathbf{u}:\Omega\times[0,1)\to\mathbb{R}^n\\
        T\in\mathbb{R}^+\\
        \mathbf{u}\mathrm{\,satisfies\,BCs\,at\,}\partial\Omega\\
        \mathbf{u}\mathrm{\,periodic\,in\,}s\mathrm{\,with\,period\,}1
    \end{array}
    \right.
    \right\}.
\end{equation}
Elements of $\mathscr{C}_p$ contain a space--time field $\mathbf{u}$ corresponding to a smooth loop parametrized by $s\in[0,1)$ in the state space, augmented by a period $T$ that reparameterizes the loop and, thereby, scales the tangent velocity vectors. The cost function $J_p^2$ takes the value of zero if and only if it is evaluated on a periodic orbit where the tangent velocity vector, $(1/T)\partial\mathbf{u}/\partial s$, coincides with the field vector, $\mathbf{f}(\mathbf{u})$, for all $s$. Following a similar logic, the identification of other types of invariant solutions, such as traveling waves \citep{farazmand2016adjoint}, relative periodic orbits \citep{parker2022variational}, invariant tori \citep{parker2023}, or connecting orbits \citep{ashtari2023connection} can be recast as a minimization problem in modified search spaces.

The landscape of the cost function continuously evolves as the bifurcation parameter is varied. In this picture, the merging and annihilation of two invariant solutions through a saddle-node bifurcation corresponds to the merging of two global minima of the cost function, that take the value of $J=0$, resulting in a minimum that lifts away from zero as the control parameter is further increased from the bifurcation value. We characterize a ghost by the minimum of the cost function with non-zero value that emerges as a result of the saddle-node bifurcation. Therefore, the ghost of an equilibrium is the locally `slowest' state, and the ghost of higher-dimensional invariant solutions are the `best-fit' sets of prescribed topological structure in the neighborhood of a destroyed solution with respect to the given cost function.

Having set up the minimization problem that yields invariant solutions and their ghosts, a suitable optimization method must be chosen from the numerous possibilities to solve the minimization problem. We employ a variational method introduced recently for computing equilibria \citep{farazmand2016adjoint, ashtari2023identifying} and periodic orbits \citep{azimi2022constructing, parker2022variational} of high-dimensional nonlinear dynamical systems. In this method, the gradient or the functional derivative of the cost function is derived as an explicit function of the unknown variables using adjoint calculations or calculus of variations. Therefore, the cost function is minimized by integrating the gradient descent dynamics
\begin{equation}
\label{eq:gradient_descent}
    \dfrac{\partial\mathbf{x}}{\partial\tau}=-\nabla_\mathbf{x} J^2,
\end{equation} 
where $\mathbf{x}$ is an element of the respective search space \eqref{eq:search_space_EQ} or \eqref{eq:search_space_PO}, and $\tau$ is a fictitious time parameterizing the gradient descent dynamics. Since $\nabla_\mathbf{x} J^2$ is derived as an explicit function of $\mathbf{x}$, the memory requirement and the computational cost of the resulting minimization method scale linearly with the size of the problem, which allows us to apply this method to high-dimensional problems including fluid flows.

This matrix-free gradient descent is not the only feasible minimization method. An alternative matrix-free implementation of the gradient descent can be achieved by using automatic differentiation \citep{page2022recurrent}. Previous authors have employed other minimization techniques to construct invariant solutions of a nonlinear dynamical system via cost function minimization, for instance, the infinitesimal-step Newton's method \citep{lan2004variational}, the Levenberg--Marquardt method \citep{parker2023} and AdaGrad \citep{page2022recurrent}. Each of these methods has advantages and disadvantages, and the choice will depend on memory limitations, robustness against inaccurate initial guesses, and the complexity of the system in question. In lower-dimensional problems, the methods based on the calculation of the Jacobian or Hessian matrices might outperform the chosen matrix-free gradient descent method in terms of speed. However, the size of the Jacobian or Hessian matrices scale quadratically with the size of the problem. This scaling is prohibitively expensive for very high-dimensional, strongly nonlinear problems such as fluid flows. Since one of our objectives is to characterize ghosts in such high-dimensional problems, we employ the aforementioned matrix-free gradient descent method throughout this paper, regardless of the size of the system studied.

The gradient descent dynamics is globally contracting, and all its trajectories eventually reach a stable equilibrium where $\partial\mathbf{x}/\partial\tau = \mathbf{0}$ and $J^2$ is minimized. Invariant solutions and their ghosts correspond to the minima of $J^2$ with zero and non-zero values, respectively; hence, they are represented by stable equilibria in the gradient descent dynamics. Consequently, the bifurcation diagram of the stable equilibria of the gradient descent dynamics incorporates ghosts into the standard saddle-node bifurcation diagram: such a bifurcation diagram consists of the folding branch corresponding to the saddle-node bifurcation, augmented by the \emph{ghost branch} that bifurcates from the saddle-node bifurcation point. The cost function is zero along the folding branch, while it takes a non-zero value along the ghost branch. It should be noted, however, that a stable equilibrium of the gradient descent dynamics is not necessarily a ghost, since a non-zero minimum of the cost function might be disconnected from any global minimum and thus from an invariant solution branch.

\section{Examples}
\label{sec:examples}
We consider a series of dynamical systems exhibiting the phenomenon of ghosts, to demonstrate the versatility of our methods. Throughout this paper, $J$ denotes the principal square root of the cost function $J^2$ defined for each problem. Note that $J$ is proportional to the (weighted) root-mean-square of the residual across all degrees of freedom in each discrete or discretized problem.

\subsection{Elastic buckling of a semi-circular arc}
\label{sec:buckling}
Elastic solid structures may suddenly jump from one equilibrium state to another when an external load exceeds a critical value---a phenomenon known as buckling in structural mechanics \cite{koiter1967stability}.
Within dynamical systems theory, many buckling problems can be interpreted as either a saddle-node or pitchfork bifurcation, with the external load being the bifurcation parameter.
However, when the load is only slightly greater than the critical value, the buckling processes can occur surprisingly slowly. 

This behavior is often termed `critical slowdown' in buckling problems \citep{gomez2017critical}. In the case of a saddle-node bifurcation, the critical slowdown can be understood as the dynamics passing the ghost of an equilibrium state.

\begin{figure}
    \centering
    \includegraphics[width=\columnwidth]{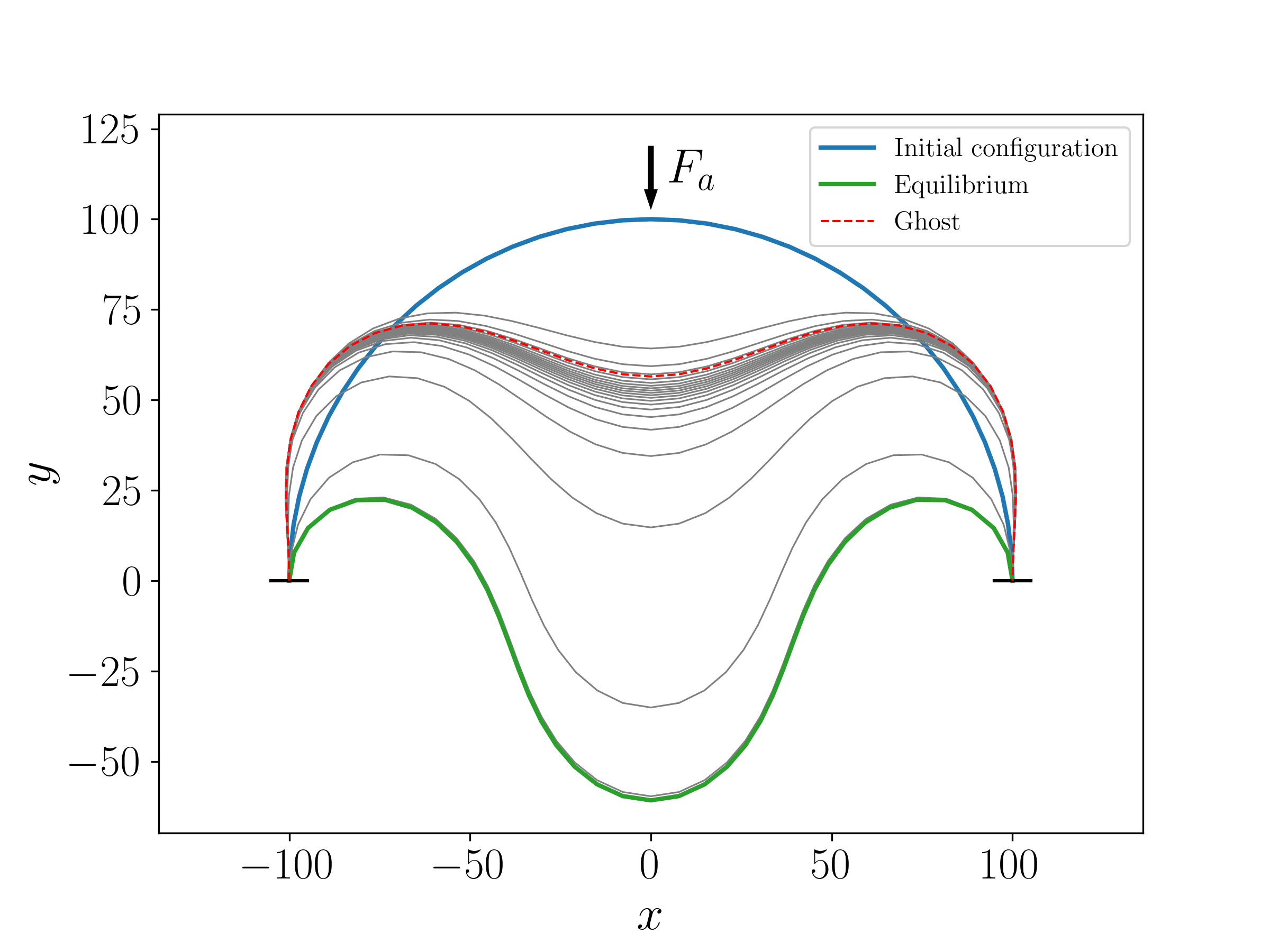}
    \caption{
    A 2D semi-circular arc with a point force $F_a$ applied vertically at the midpoint, while the two ends of the arc are clamped.
    Solid blue and green curves indicate the initial and final configurations, respectively. The gray curves demonstrate the evolution of the arc at uniform time intervals of $\Delta t = 2000$ for $F_a = 1060$, that is slightly above the buckling load.
    Notice the slowdown near the ghost state, marked by a dashed red curve.
    }
    \label{fig:arc_shape}
\end{figure}

Here, we consider a two-dimensional (2D) semi-circular arc whose buckling is of the saddle-node bifurcation type.
The two ends of the arc are clamped, and a point force $F_a$ is applied at the midpoint, directed toward the center of the arc’s defining circle (see Fig.~\ref{fig:arc_shape}). We describe the arc using 2D Reissner beam theory \citep{reissner1972one} and derive the nonlinear equations of motion from Hamilton's principle \citep{lai2009introduction, ibrahimbegovic1999nonlinear}.

Following a standard finite element method (FEM), the displacements in the $x$ and $y$ directions, as well as the rotation of the arc, are interpolated using Lagrange polynomials. As a result, we obtain the following semi-discrete equation of motion in matrix form:
\begin{equation}
    \mathbf{M} \ddot{\mathbf{U}}+ \mathbf{D} \dot{\mathbf{U}} + \mathbf{F}^\text{int}(\mathbf{U})  = \mathbf{F}^{\text{ext}}, \label{eq:arc_motion}
\end{equation}
where $\mathbf{M}$ and $\mathbf{D}$ are the constant mass and damping matrices, respectively, which depend on the geometry and material properties of the arc; $\mathbf{U}$ denotes the generalized displacement vector of all nodes; the vector $\mathbf{F}^{\text{int}}$ represents the internal force at the nodes, which is a nonlinear function of  $\mathbf{U}$; and $\mathbf{F}^{\text{ext}}$ is the external force vector at the nodes.

We study the dynamics of buckling of the semi-circular arc based on the above beam model using 40 linear elements. The material and geometric parameters are the arc radius $R = 100$, Young's modulus $E = 1 \times 10 ^6$, second moment of inertia $I = 1$, cross-sectional area $A = 2.29$, Poisson's ratio $\nu = 0$, density $\rho = 50$, and damping coefficient $\mu = 1$, all in consistent physical units. These values are taken from the verification example in Chap. 13 of the book \cite{zienkiewicz2014finite} and do not represent any real material or structure.
For time integration of Eq.~\eqref{eq:arc_motion}, we employ the implicit Newmark-$\beta$ method, which ensures numerical stability (see Ref.~\citep{simo1986dynamics} for details).
In a numerical simulation for a given value of the external force, $F_a$ increases linearly from zero to the target value over a time interval of $t = 100$, after which it remains constant.

\begin{figure}
    \centering
    \includegraphics[width=1\columnwidth]{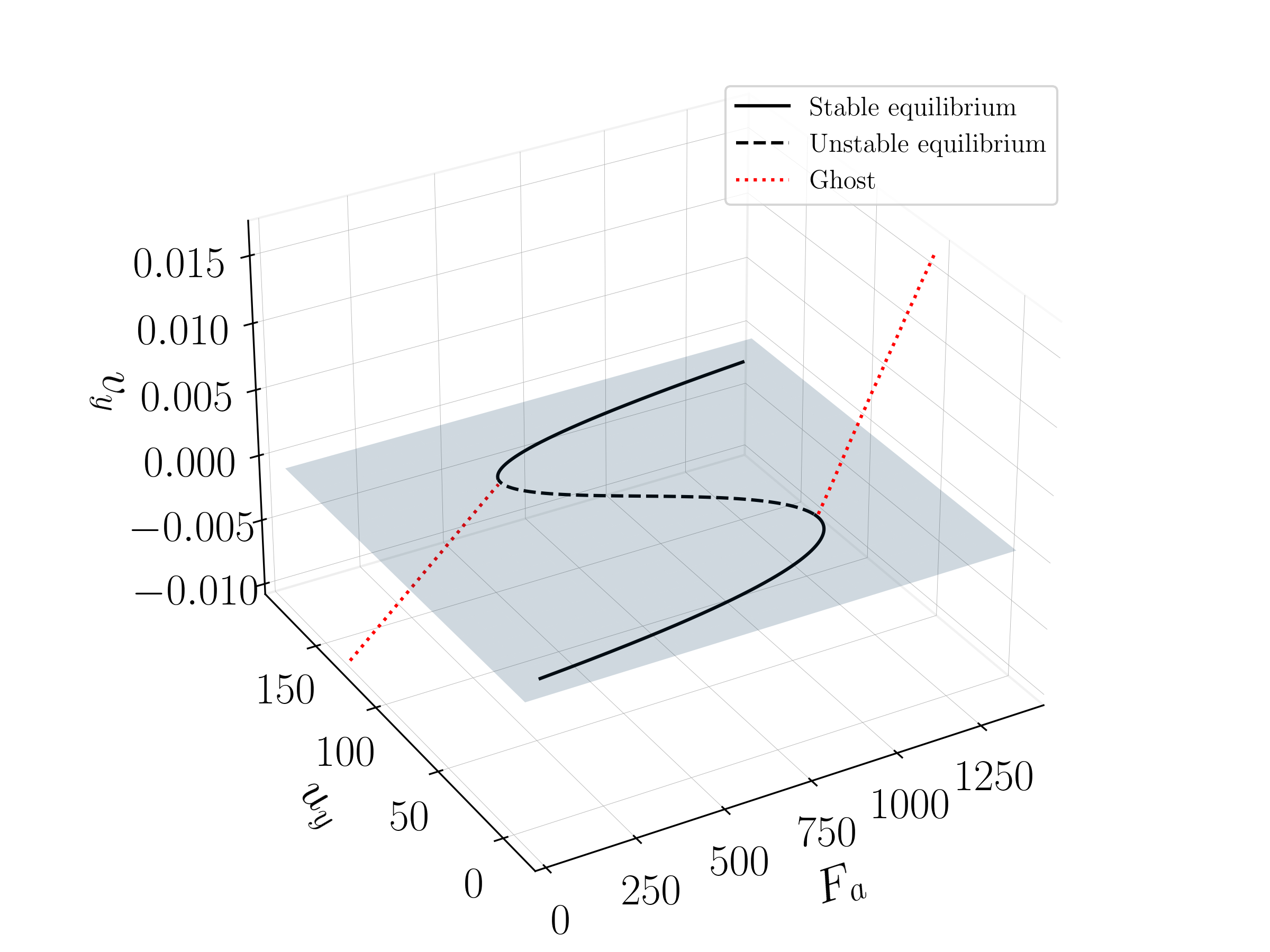}
    \caption{
    Bifurcation diagram for the 2D semi-circular arc system. Stable and unstable equilibria are represented by solid and dashed black lines, respectively, while the ghost branches, continued from the two saddle-node bifurcation points, are shown as red dotted lines.
    The vertical displacement $u_y$ and vertical velocity $v_y$ of the arc's midpoint are both plotted as functions of the external force $F_a$.
    The equilibrium solution branch lies within the zero-velocity plane, while the ghost branch extends outside this plane, following the slowest state with respect to the cost function \eqref{eq:cost_arc}.}
    \label{fig:arc_ghost}
\end{figure}

For any nonzero damping force, the dynamics can be expressed in variational form. As energy is dissipated, the arc dynamics must converge to an equilibrium state for a fixed $F_a$; in other words, periodic solutions are not allowed.
The possible vertical displacements $u_y$ of the arc's midpoint at equilibrium, as a function of the external force $F_a$, are summarized in Fig.~\ref{fig:arc_ghost}. Two saddle-node bifurcations exist for $F_a = 445.4$ and $F_a = 1056.8$.
The larger value is called the `buckling force,' as for larger forces, the only possible equilibrium state is on the upper branch; hence, a transition from the lower to the upper branch takes place once the external force exceeds this critical value.

\begin{figure}
    \centering
    \includegraphics[width=1\columnwidth]{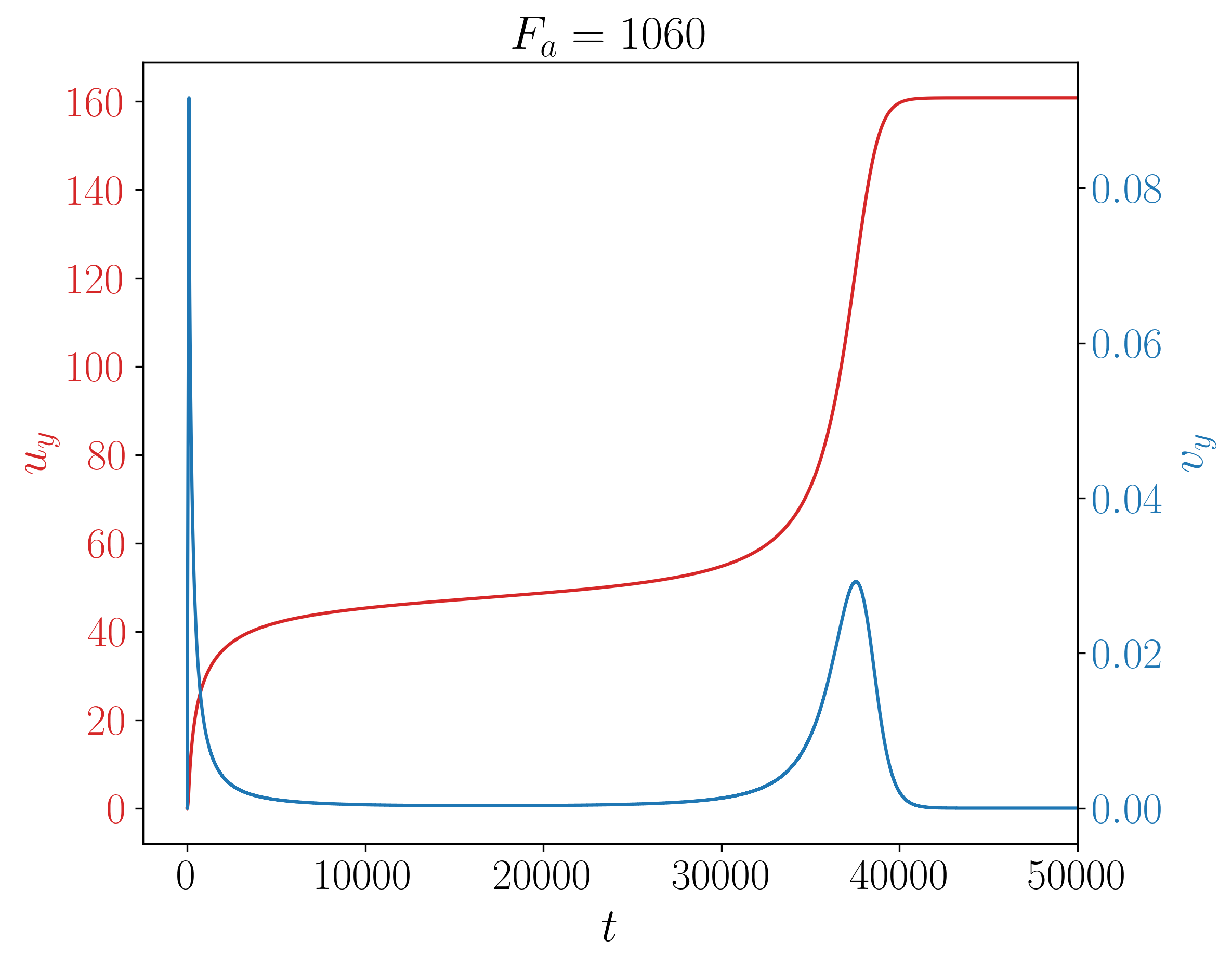}
    \includegraphics[width=1\columnwidth]{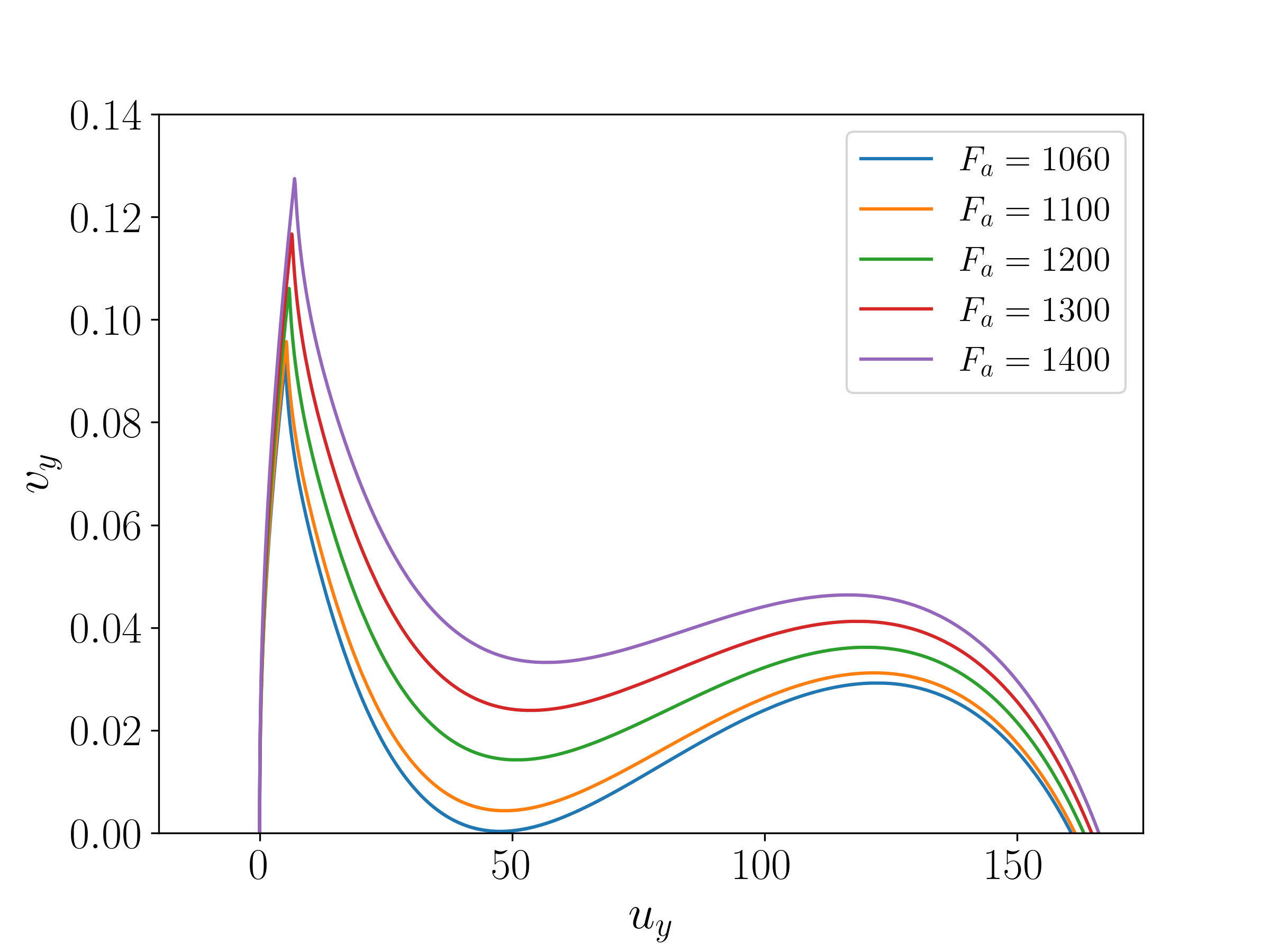}
    \caption{Elastic buckling of the 2D semi-circular arc for different values of the external load $F_a$.
    Above: Time evolution of the vertical displacement $u_y$ and vertical velocity $v_y$ of the arc's midpoint for $F_a=1060$, which is slightly above the saddle-node bifurcation value.
    An extended period of near-zero velocity is a signature of the ghost phenomenon.
    Below: The trajectories of the system projected onto the velocity--displacement plane of the arc's midpoint for different values of $F_a$. The slowdown due to the ghost is evident, even for values of $F_a$ much larger than the bifurcation value.}
    \label{fig:arc_DNS}
\end{figure}

Figure \ref{fig:arc_DNS} shows simulation results for various values of the external force. For $F_a = 1060$, which is slightly above the bifurcation value, we observe a plateau in the time evolution of the vertical displacement $u_y$, corresponding to an extended period of near-zero vertical velocity $v_y$ (see the top panel of Fig.~\ref{fig:arc_DNS}). These are signatures of the ghost phenomenon or critical slowdown.
A clear view of this phenomenon can also be seen in the evolution of the arc's configuration during this simulation, as shown in Fig.~\ref{fig:arc_shape}. In this figure, the density of intermediate snapshots, uniformly spaced in time, reveals a bottleneck through which the system takes a substantial amount of time to pass.
The slowed transition due to the ghost of the bifurcation remains observable even for values of $F_a$ considerably larger than the buckling force, as evident from the bottom panel of Fig.~\ref{fig:arc_DNS}.

In the following, the representative ghost state is defined and computed for this system.

\subsubsection{Definition and computation of the ghost states}
In order to define the ghost state following the methodology introduced in Sec. \ref{sec:method}, we rewrite Eq.~\eqref{eq:arc_motion} as the following set of first-order evolution equations:
\begin{subequations}
    \begin{align}
        \dot{\mathbf{U}} & = \mathbf{V}, \\
        \dot{\mathbf{V}} & = \mathbf{M}^{-1}(\mathbf{F}^\text{ext} - \mathbf{F}^\text{int}  - \mathbf{D} \mathbf{V}).
    \end{align}
    \label{eq:first_order_motion}
\end{subequations}
This equation has the form of the general dynamical system \eqref{eq:general_governing_equation} with $d=0$ and $n= 6n_\text{node}$, where $n_\text{node}$ is the number of nodes in the FEM discretization. The state at each node is determined by six degrees of freedom: two displacements in the $x$ and $y$ directions, one rotation in the $x$--$y$ plane, and the three corresponding velocities. According to the method described in Sec. \ref{sec:method}, we define the cost function $J$ for fixed points and derive its gradient descent dynamics as follows:
\begin{equation}
\label{eq:cost_arc}
     J^2 := \frac{1}{2} \left(\mathbf{\dot{U}}^\text{T} \mathbf{\dot{U}} + \mathbf{\dot{V}}^\text{T} \mathbf{\dot{V}}\right), 
 \end{equation}
\begin{equation}
     \frac{\mathrm{d}}{\mathrm{d} \tau} \begin{bmatrix} 
     \mathbf{U} \\
     \mathbf{V}
     \end{bmatrix} 
     = - \nabla J^2  
    = -\begin{bmatrix}
        - \mathbf{K} \mathbf{M}^{-1} \dot{\mathbf{V}} \\ 
        \mathbf{V} - \mu \dot{\mathbf{V}}
    \end{bmatrix},
    \label{eq:arc_gradient_descent}
\end{equation}
where $\mathbf{K} = {\partial \mathbf{F}^\text{int}}/{\partial \mathbf{U}}$ is the tangent stiffness matrix. 

The global minima of the cost function $J$, which take the value of zero, correspond to the equilibrium states, while the nonzero minima that appear following the saddle-node bifurcations correspond to the ghost states.
In order to minimize $J$, we integrate the gradient descent dynamics \eqref{eq:arc_gradient_descent} forward in the fictitious time $\tau$. If the minimum value of $J^2$ falls below $10^{-20}$, we consider the obtained state as a fixed point.

Starting from the equilibrium state at each of the two saddle-node bifurcations, we track the ghost states by increasing the external force from the buckling load $F_a=1056.8$ or by decreasing it from $F_a=445.4$. For the ghost states, we observe that the displacement $u_y$ of the arc's midpoint remains almost constant across different values of the external force (see Fig.~\ref{fig:arc_ghost}).
This is because the arc's nonlinear internal force is locally extremal for that displacement. The extremal internal force is a maximum for the ghosts of the bifurcation at $F_a=1056.8$, and a minimum for the ghosts of the bifurcation at $F_a=445.4$. In either case, the arc's internal force counteracts and resists the external force; therefore, the minimum acceleration is achieved around the particular configuration of the ghost states, regardless of the chosen value of $F_a$. However, the velocity $v_y$ of the arc's midpoint increases as we follow the ghost states further away from the bifurcation point (see Fig.~\ref{fig:arc_ghost}). Correspondingly, the minimum value of $J$ increases along the ghost branch as the control parameter is varied further away from the bifurcation point. 

We conclude that critical slowdown in buckling problems can be interpreted as the ghost of equilibrium states, remaining in the state space after their destruction through a saddle-node bifurcation. By defining and minimizing a suitable non-negative cost function, we are able to compute both equilibrium and ghost states, enabling us to characterize this phenomenon quantitatively in buckling problems.

\subsection{One-dimensional polynomial map}
\label{sec:map}
\begin{figure}
    \centering
    \includegraphics[width = 0.8\columnwidth]{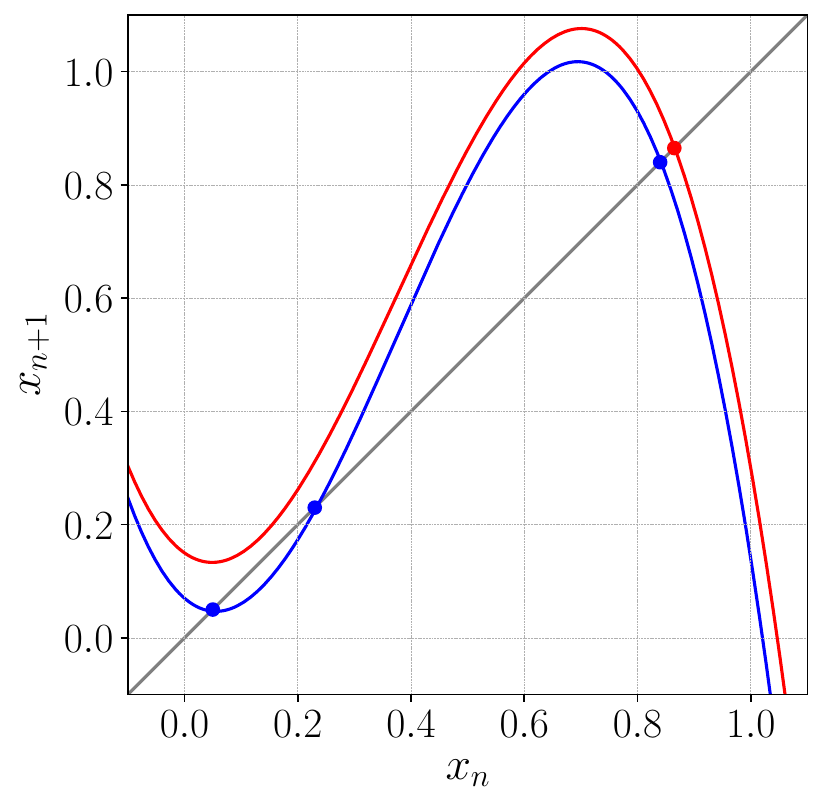}
    \caption{The map (\ref{eq:map}) at parameter values $\rho=0.07$ (blue) and $\rho=0.15$ (red). In the former case, below the bifurcation value, there are three fixed points. In the latter case, above the bifurcation value, there is only one.}
    \label{fig:map}
\end{figure}

As an example of a chaotic system, consider the map defined by the recurrence relation
\begin{equation}
x_{n+1}=f(x_n):=(\rho-1)(8x_n^3-9x_n^2)+2\rho x_n - x_n + \rho.
\label{eq:map}
\end{equation}
Here, $\rho$ is the bifurcation parameter. At $\rho^\ast = \left(33\sqrt{33}-115\right)/668$, the system undergoes a saddle-node bifurcation; below this value, there are three fixed points, and immediately above it, there is only one (see Fig.~\ref{fig:map}).

\begin{figure}
    \centering
    \includegraphics[width=\columnwidth]{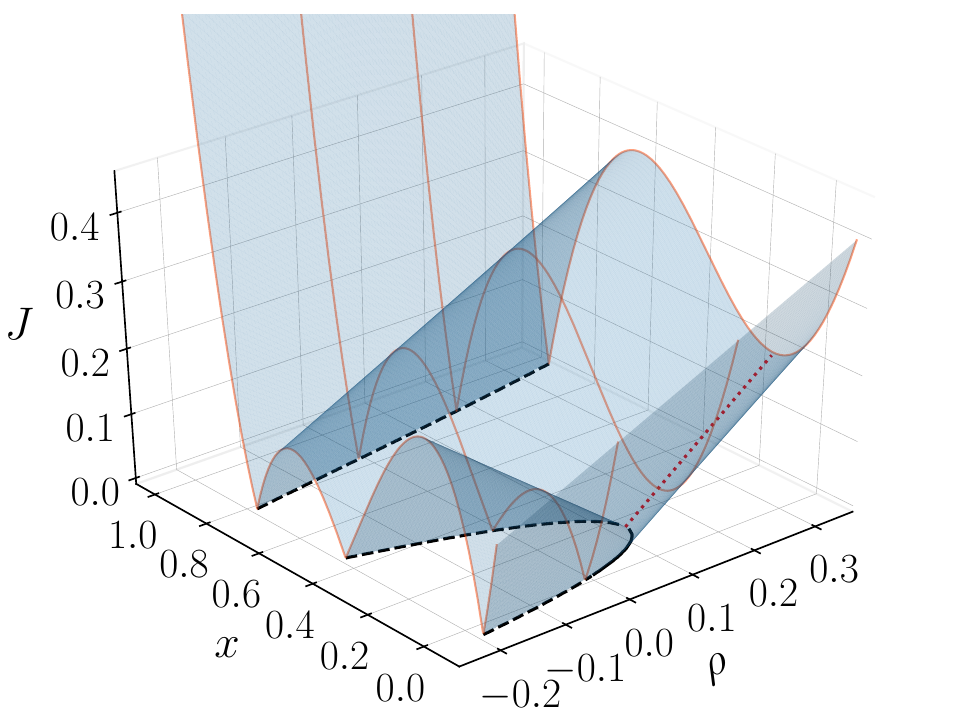}
    \caption{The landscape of the cost function $J$ (Eq.~\eqref{eq:mapJ}) for the polynomial map \eqref{eq:map}. The global minima correspond to fixed points, represented by solid lines for the stable and dashed lines for the unstable ones. The dotted line marks the non-zero local minima that correspond to the ghost of the saddle-node bifurcation. Above the bifurcation point, $J$ increases linearly with $\rho$.}
    \label{fig:extrema}
\end{figure}

This system is simple enough to solve for fixed points analytically without defining and minimizing a cost function, as discussed in Sec. \ref{sec:method}. However, we choose to do so anyway to demonstrate the method. A cost function that penalizes the deviation of $x$ from being a fixed point of the map, i.e., the discrete-time equivalent of Eq.~\eqref{eq:ghosts_cost_FP}, is
\begin{equation}
    J^2(x) := \frac{1}{2}\left(f\left(x\right)-x\right)^2.
    \label{eq:mapJ}
\end{equation}
This has gradient ${\mathrm{d}J^2}/{\mathrm{d}x} = \left(f(x)-x\right)\left(f'(x)-1\right)$. Hence, extrema exist at fixed points $f(x)=x$, as expected, but also at $f'(x)=1$. From Eq.~\eqref{eq:map},
\begin{equation}
    f'(x)-1 = (\rho-1)(24x^2-18x+2),
\end{equation}
which has roots at $(9\pm\sqrt{33})/24\approx0.614,0.136$. These extrema of the cost function are illustrated in Fig.~\ref{fig:extrema}. We see that the larger of these roots is always a maximum for the range of parameters considered, but the smaller one changes from a maximum to a minimum after the saddle-node bifurcation.

\begin{figure}
    \centering
    \includegraphics[width=\columnwidth]{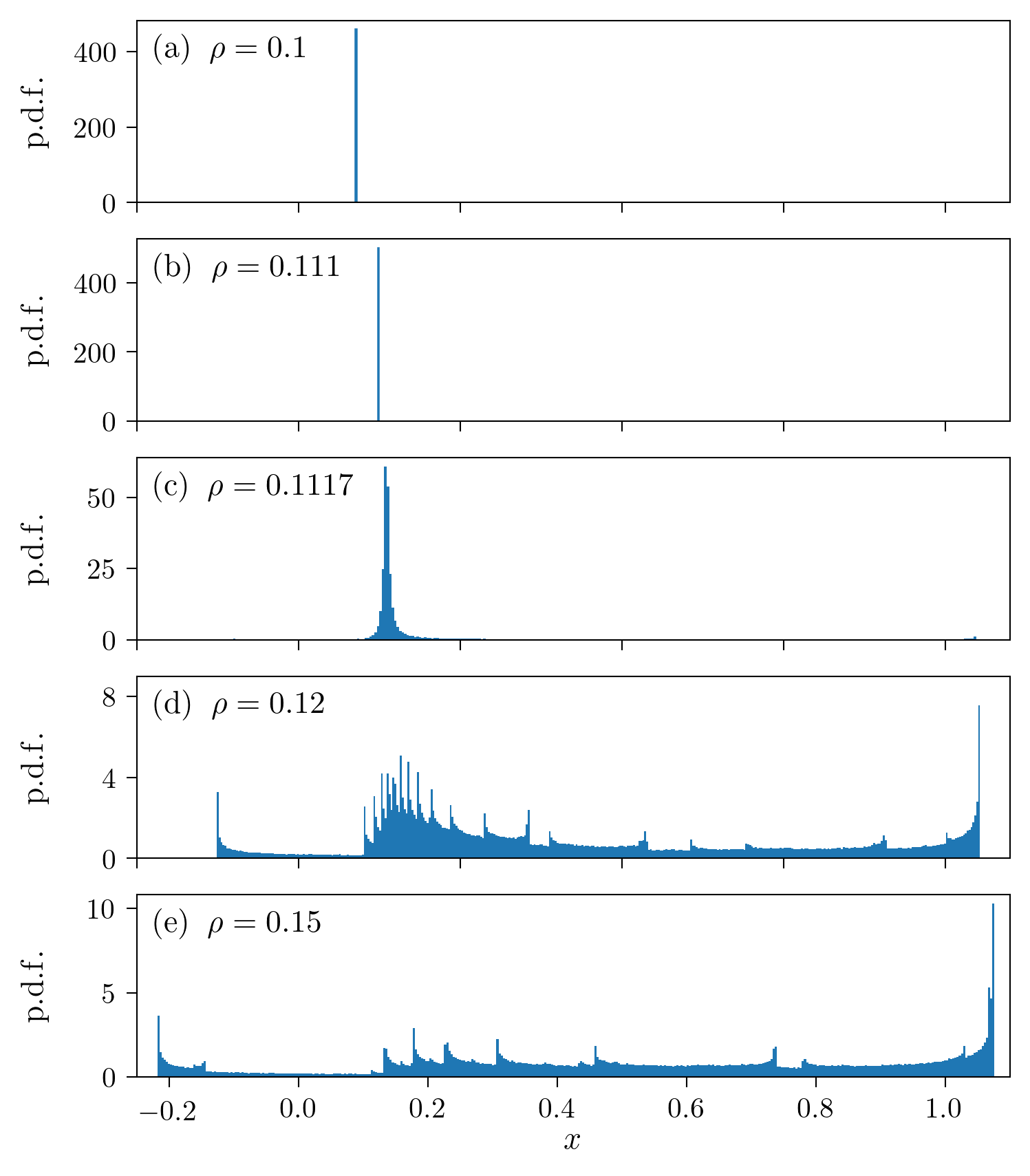}
    \caption{Monte Carlo approximations of probability density functions (pdfs) of the solution to the 1D polynomial map (\ref{eq:map}). Above the bifurcation, the map is chaotic, and the pdfs are fractal. (a,b) Below the bifurcation, for $\rho=0.1$ and $\rho=0.111$, the long-time solution of the system is exactly the attracting fixed point. Therefore, the pdf is a delta function at its location. (c) Just above bifurcation, for $\rho=0.1117$, the ghost creates a peak in the pdf, and the trajectory spends nearly every iteration very close to the minimum of the cost function, i.e., the ghost state. (d) Further from the bifurcation, for $\rho=0.12$, the effect of the ghost on the dynamics is seen as a peak in the pdf. (e) Further still from the bifurcation, for $\rho=0.15$, the ghost no longer has a discernible effect on the dynamics.}
    \label{fig:1dpdfs}
\end{figure}

The long-term behavior of the map for different values of $\rho$ is visualized in Fig.~\ref{fig:1dpdfs} in terms of the probability density function (pdf) of a sufficiently long sequence generated by the map. Below the bifurcation, for $0<\rho<\rho^\ast$, trajectories starting sufficiently close are attracted to the one stable fixed point. Therefore, the pdf becomes a delta function at the location of the stable solution (see panels (a) and (b)). Above the bifurcation, the map is chaotic, and the pdf densely fills the plotted range. However, this chaos is intermittent, and for parameter values close to the bifurcation value, trajectories spend many iterations in the vicinity of the local minimum with a non-zero value. Consequently, the pdf exhibits a peak at the local minimum of the cost function, i.e., the ghost state (see panel (c)). Note that this is not true for all minima: the unstable fixed points also manifest as minima of the cost function, but no trace of them at all is seen in the chaotic pdf. As $\rho$ is further increased, the ghost gets weaker (see panel (d)) until, sufficiently far from the bifurcation point, it is no longer able to create a pronounced peak in the pdf (see panel (e)).

In this simple example, we can also consider the complex roots to $f(x)=x$ as a different approach for studying ghosts (see Sec. \ref{sec:intro}). In this case, above the bifurcation, we find two additional complex solutions. For example, for $\rho=0.12$,
\begin{equation}
    x\approx0.134355 \pm 0.0430683 i.
\end{equation}
The real part of this expression could be associated with the ghost, following \citet{canela2022dynamical}. Observe that now the value of this real part depends on the value of the parameter $\rho$, which was not the case for the method discussed above. However, close to the bifurcation, where the effect of the ghost is most strongly felt, the positions are very similar between these different definitions.

\subsection{The Lorenz system}
\label{sec:lorenz}
\begin{figure}
    \includegraphics[width=\columnwidth]{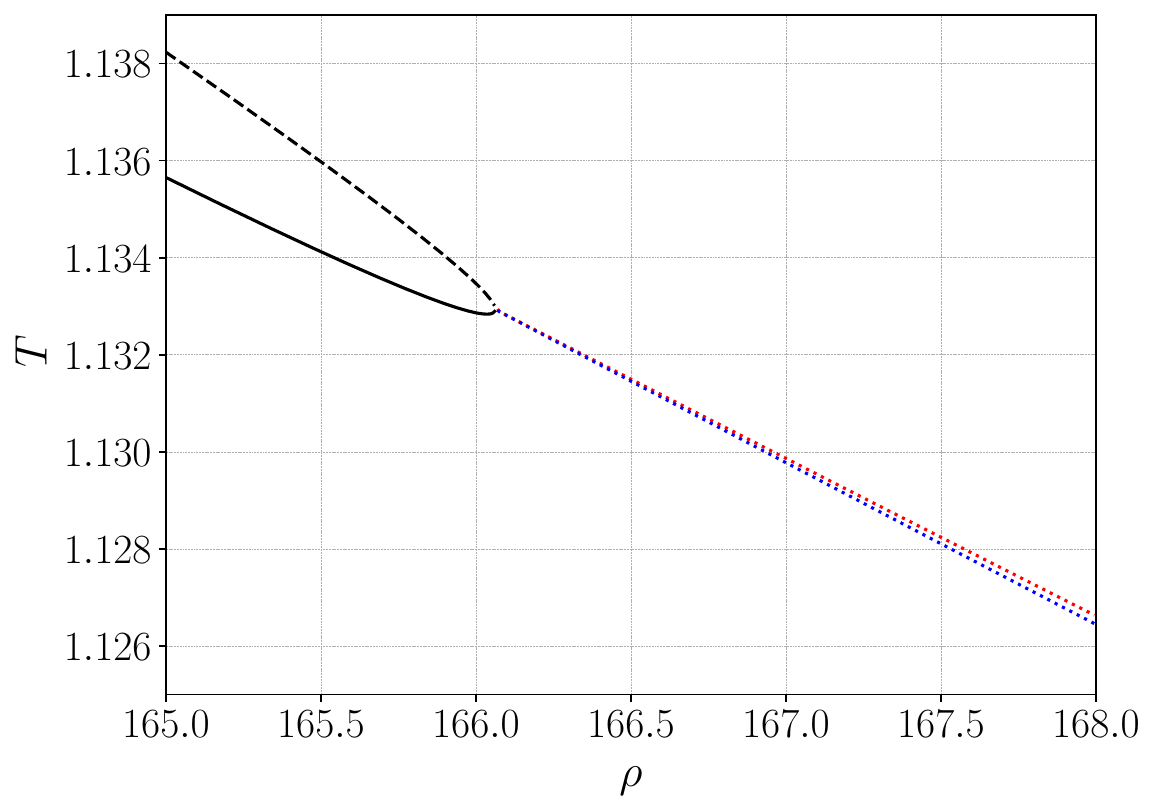}
    \centering
    \caption{Bifurcation diagram for the Lorenz system, with $\sigma=10$, $\beta=8/3$, and $\rho$ varying. Below $\rho\approx166.052$, a stable and an unstable periodic orbit exist, marked by solid and dashed black lines, respectively. These annihilate in a saddle-node bifurcation, resulting in a ghost state, which is marked by the dotted lines. The vertical axis shows the period of the structures. The blue dotted line shows the ghost calculated by using the cost function $J_1$, and the red dotted line instead uses $J_{100}$ (see the definition \eqref{eq:ghosts_J_alpha}).}
    \label{fig:lorenzbifurc}
\end{figure}

As a simple system of ODEs with a chaotic attractor, we consider Lorenz's famous 1963 system \citep{lorenz1963deterministic}:
\begin{equation}
\dfrac{\mathrm{d}}{\mathrm{d}t}
    \begin{bmatrix}
        x\\
        y\\
        z
    \end{bmatrix} = 
    \begin{bmatrix}
        \sigma(y-x)\\
        x(\rho-z)-y\\
        xy-\beta z
    \end{bmatrix}.
    \label{eq:lorenz}
\end{equation}
Here, $\sigma$, $\beta$, and $\rho$ are constant parameters of the system, originally set by Lorenz to $10$, $8/3$, and $28$, respectively. Instead of Lorenz's original settings, \citet{manneville1979intermittency} performed numerical simulations with $\sigma=10$, $\beta=8/3$, and $\rho\approx 166$. At $\rho=166$, they observed perfectly periodic behavior. However, at $\rho=166.1$, this same periodic pattern is interrupted, seemingly at random, by incoherent chaotic behavior. This is an example of intermittent chaos, similar to that in the previous example, but now in a continuous-time dynamical system. The chaos appears to intermittently visit a periodic orbit; however, upon continuation of the attracting orbit at $\rho=166$, we observe that it undergoes a saddle-node bifurcation at $\rho\approx166.052$ (see Fig.~\ref{fig:lorenzbifurc}). This indicates that there is no relevant invariant solution at $\rho=166.1$; instead, intermittent chaos is shaped by the ghost of the saddle-node bifurcation.

\begin{figure*}
    \centering
    \includegraphics[width=\textwidth]{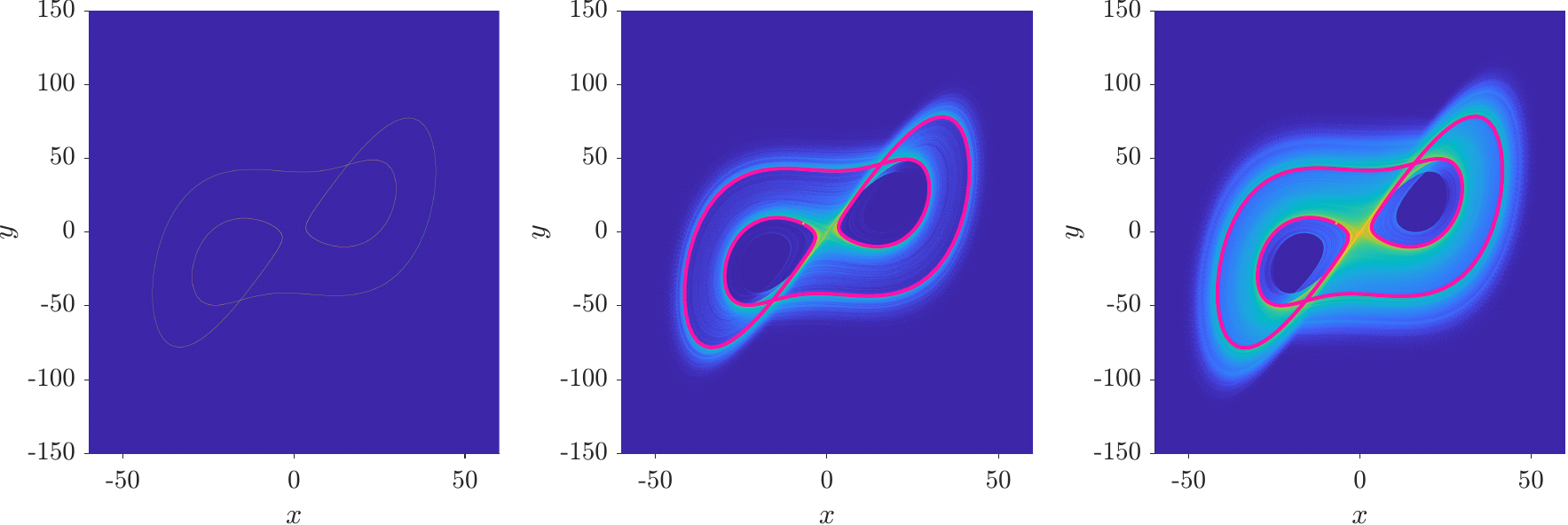}
    \includegraphics[width=\textwidth]{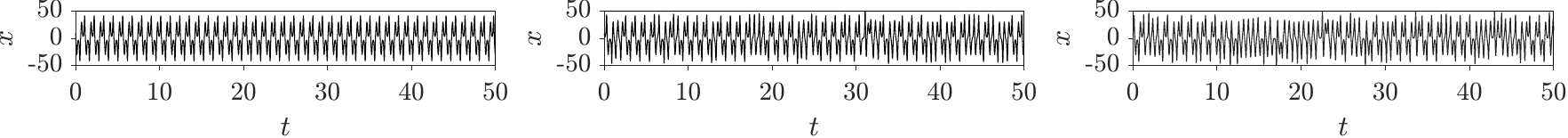}
    \caption{Pdfs and typical time series for the attractor in the Lorenz system with $\sigma=10$, $\beta=8/3$, and different values of $\rho$. Left: $\rho=166$, below the bifurcation; the entire pdf collapses onto the 1D stable periodic orbit. The resulting line of non-zero probability is thus barely visible. Middle: $\rho=167$, above the bifurcation; the attractor is chaotic, but it is focused around the ghost of the periodic orbit, shown in red. Right: $\rho=168$; the chaotic attractor is more complicated, but the ghost is still visible (see Fig.~\ref{fig:lorenzzoom}).}
    \label{fig:lorenzresults}
\end{figure*}

\begin{figure*}
    \centering
    \includegraphics[width=0.9\columnwidth]{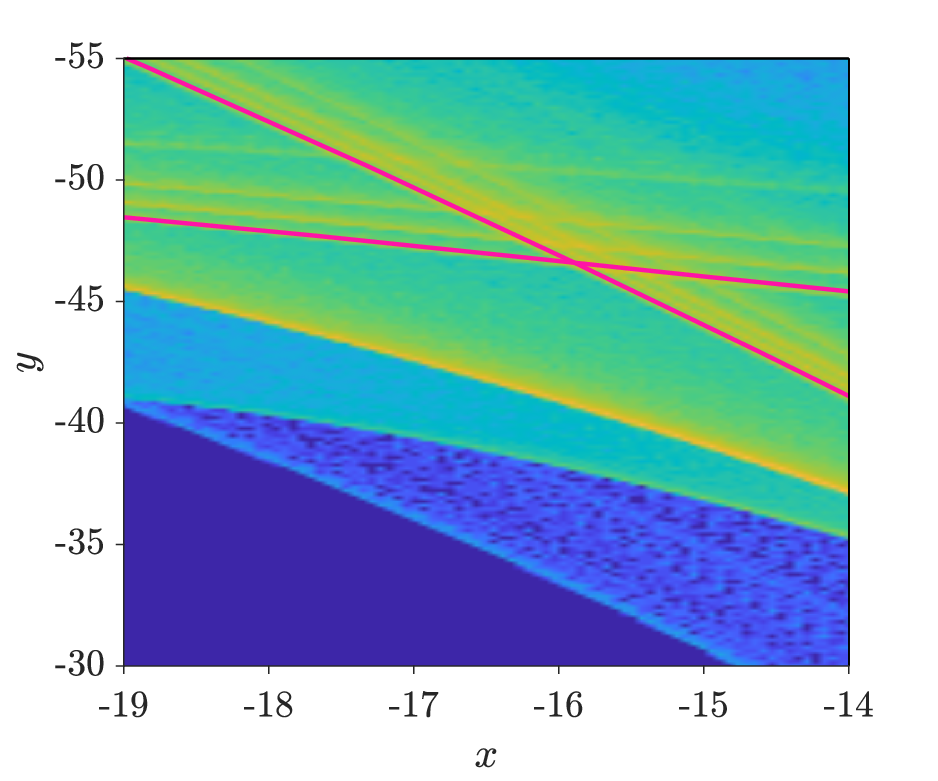}
    \includegraphics[width=0.9\columnwidth]{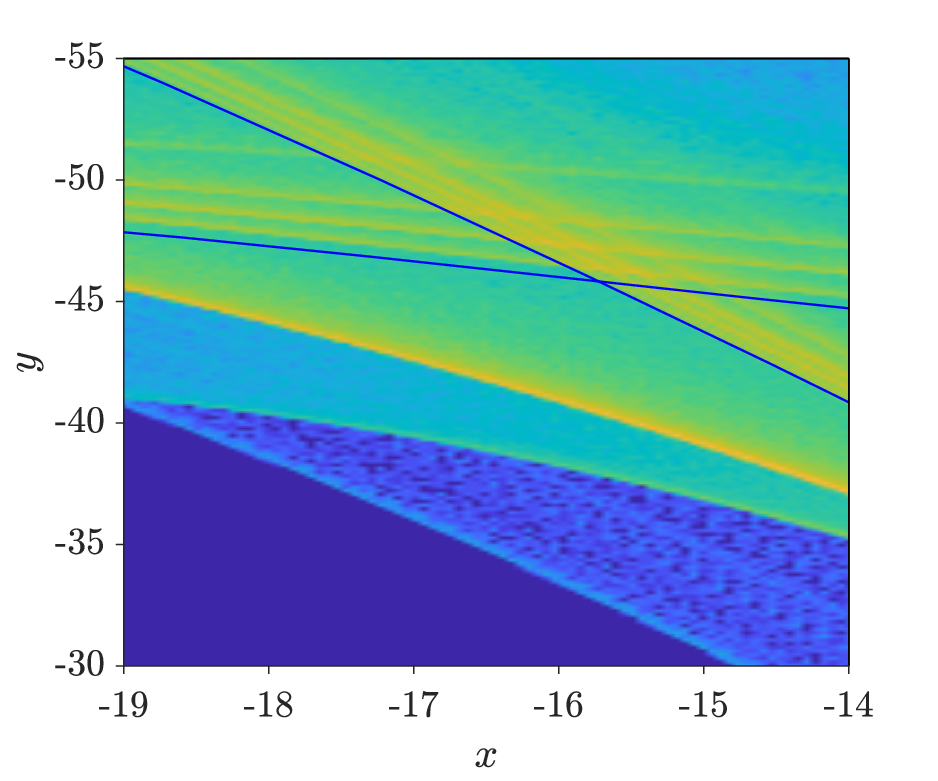}
    \caption{Zoom of Fig.~\ref{fig:lorenzresults} showing the pdf in the state space of the Lorenz dynamics for $\rho=168$. Left: The ghost found by our method, shown in red, obstructs one of the frequently visited paths, represented by yellow-hued stripes, indicating a close match.
    Right: The stable periodic orbit at $\rho=166$, overlaid in blue, does not closely match any of the yellow stripes, demonstrating that the continuation of the ghost is not merely the periodic orbit itself. While neither the periodic orbit at a different parameter value nor its ghost continued to the present parameter value is an invariant solution of the dynamics at $\rho=168$, the latter matches the statistical structure of the attractor more precisely.}
    \label{fig:lorenzzoom}
\end{figure*}

The Lorenz system has the form of the general dynamical system \eqref{eq:general_governing_equation} with $d=0$ and $n=3$. Following the methods discussed in Sec. \ref{sec:method}, an obvious choice of cost function to find periodic orbits in this system is
\begin{multline}
J_1^2\left[x(\cdot),y(\cdot),z(\cdot),T\right]:= \\
\begin{aligned}
&\frac{1}{2}\int_0^{1} \left(\frac{1}{T}\frac{\mathrm{d}x}{\mathrm{d}s} - \sigma\left(y-x\right)\right)^2 \mathrm{d}s  \\
+&\frac{1}{2}\int_0^{1} \left(\frac{1}{T}\frac{\mathrm{d}y}{\mathrm{d}s} - x\left(\rho-z\right) + y\right)^2 \mathrm{d}s \\
+&\frac{1}{2}\int_0^{1} \left(\frac{1}{T}\frac{\mathrm{d}z}{\mathrm{d}s} - xy+\beta z\right)^2 \mathrm{d}s.
\end{aligned}
\label{eq:J1}
\end{multline}
However, there are many other possible choices of cost function, for example
\begin{multline}
J_{\alpha}^2\left[x(\cdot),y(\cdot),z(\cdot),T\right] := \\
\begin{aligned}
&\frac{1}{2}\int_0^{1} \left(\frac{1}{T}\frac{\mathrm{d}x}{\mathrm{d}s} - \sigma\left(y-x\right)\right)^2 \mathrm{d}s \\
+&\frac{1}{2}\int_0^{1} \left(\frac{1}{T}\frac{\mathrm{d}y}{\mathrm{d}s} - x\left(\rho-z\right) + y\right)^2 \mathrm{d}s \\
+&\frac{1}{2\alpha^2}\int_0^{1} \left(\frac{1}{T}\frac{\mathrm{d}z}{\mathrm{d}s} - xy+\beta z\right)^2 \mathrm{d}s,
\label{eq:ghosts_J_alpha}
\end{aligned}
\end{multline}
for a real constant $\alpha\neq 0$. Such a cost function also describes a periodic orbit precisely when it takes the global minimum value of zero. However, the exact properties of the ghost can change, as the cost function is minimized at a different loop in state space for different choices of $\alpha$. Figure \ref{fig:lorenzbifurc} shows the ghost in the bifurcation diagram for two different choices of $\alpha$. The usual choice, $\alpha=1$, represents a balance between the extremes $\alpha\to0$ (not pictured) and $\alpha\to\infty$ (well represented by $\alpha=100$). At the bifurcation point in the parameter space, the ghosts from different possible cost functions must agree, but as $\rho$ is increased, they can diverge. More exotic norms, such as using a fourth power instead of a square in Eq.~\eqref{eq:J1}, were not found to make a significant difference to the results in this particular case. They do, however, affect the performance of the computations. \citet{parker2022variational} compared different possible choices of cost function for converging periodic orbits in the 2D Navier--Stokes equations, and found some choices to be much more efficient than others.

Using $J_1$ as the cost function, the gradient-descent dynamics for the evolution of a loop, parametrized by $s\in[0,1)$, is derived as
\begin{multline}
    \dfrac{\partial}{\partial \tau}
    \begin{bmatrix}
        x(s;\tau)\\
        y(s;\tau)\\
        z(s;\tau)\\
        T(\tau)
    \end{bmatrix} = \\
    \renewcommand\arraystretch{2.1}
    \begin{bmatrix}
        \dfrac{1}{T}\dfrac{\partial r_x}{\partial s} - \sigma r_x + (\rho-z) r_y + yr_z\\
        \dfrac{1}{T}\dfrac{\partial r_y}{\partial s} + \sigma r_x - r_y + xr_z\\
        \dfrac{1}{T}\dfrac{\partial r_z}{\partial s} - xr_y - \beta r_z\\
        \dfrac{1}{T^2}\displaystyle\int_0^1\left(\dfrac{\partial x}{\partial s}r_x + \dfrac{\partial y}{\partial s}r_y + \dfrac{\partial z}{\partial s}r_z \right)\,\mathrm{d}s
    \end{bmatrix},
    \renewcommand\arraystretch{1}
    \label{eq:ghosts_Lorenz_adjoint}
\end{multline}
where $[r_x, r_y, r_z]^\top$ is the residual of the Lorenz dynamics:
\begin{equation}
    \begin{bmatrix}
        r_x\\
        r_y\\
        r_z
    \end{bmatrix} :=
    \dfrac{1}{T}\dfrac{\partial}{\partial s}
    \begin{bmatrix}
        x\\
        y\\
        z
    \end{bmatrix} - 
    \begin{bmatrix}
        \sigma(y-x)\\
        x(\rho-z)-y\\
        xy-\beta z
    \end{bmatrix}.
\end{equation}
The derivation follows the method discussed in Ref.~\cite{azimi2022constructing}. The resulting ghost as $\tau\to\infty$ is depicted in Fig.~\ref{fig:lorenzresults}. Below the saddle-node bifurcation, when the stable periodic orbit exists, all trajectories converge to it. Above the bifurcation, the fractal pdf indicates that while the dynamics is chaotic, it regularly visits the ghost. Furthermore, we see that tracking the ghost is essential: merely taking the periodic orbit at the bifurcation point provides a worse match for the highlighted region in the pdf. This is demonstrated in Fig.~\ref{fig:lorenzzoom}.

\subsection{The Kuramoto--Sivashinsky equation}
While the examples in the previous sections involved ODEs, we now consider a nonlinear PDE and study the ghosts of its periodic orbits. Specifically, we consider the 1D Kuramoto--Sivashinsky equation (KSE) \cite{Kuramoto1976,sivashinsky1977flames}. The KSE is relevant in several physical contexts, including the dynamics of flame fronts \cite{sivashinsky1977flames}, plasma physics \cite{LaQuey1975}, and interfacial fluid instability \cite{Hooper1985}, among others. This equation is believed to be the simplest PDE that exhibits spatiotemporal chaos \cite{brummitt2009search}. Therefore, it is often used as a sand-box model to test concepts related to and methods developed for nonlinear and chaotic dynamical systems \citep{Christiansen1997, cvitanovic2010state, lasagna2018sensitivity, Pathak2018, Zeng2022}.

The 1D KSE can be written as
\begin{equation}
\label{eq:KSE}
    \dfrac{\partial u}{\partial t} = - u\dfrac{\partial u}{\partial x} - \dfrac{\partial^2 u}{\partial x^2} - \dfrac{\partial^4 u}{\partial x^4},
\end{equation}
for a real-valued field $u(x,t)$, subject to periodic BCs in the spatial dimension $x$ with period $L$. In this formulation, the domain length $L$ is the only control parameter of the system. The KSE is equivariant under reflection about the origin,
\begin{equation}
    u(x,t)\mapsto-u(-x,t).
\end{equation}
We set the domain length to $L=39$ and study the KSE dynamics in the anti-symmetric subspace of functions such that $u(x,t)=-u(-x,t)$. The imposed discrete symmetry significantly reduces the complexity of the dynamics \citep{edson2019lyapunov}, so that low-dimensional chaos is observed for the chosen domain length.

The KSE \eqref{eq:KSE} has the form of the general autonomous dynamical system \eqref{eq:general_governing_equation} with $n=d=1$.
We define the relevant cost function for a loop $[u(x,s),T]\in\mathscr{C}_p$ as (see Sec. \ref{sec:method})
\begin{equation}
    J^2 := \dfrac{1}{2}\int_0^1\int_0^L r^2\,\mathrm{d}x\mathrm{d}s,
    \label{eq:KSE_J}
\end{equation}
where $r$ is the residual of the KSE:
\begin{equation}
    r := \dfrac{1}{T}\frac{\partial u}{\partial s} - \left(- u\dfrac{\partial u}{\partial x} - \dfrac{\partial^2 u}{\partial x^2} - \dfrac{\partial^4 u}{\partial x^4}\right).
\end{equation}
The gradient-descent dynamics of the cost function, governing the evolution of a loop, is given by (see Ref.~\citep{azimi2022constructing} for the derivation)
\begin{multline}
    \label{eq:KSE_loop_dynamics}
    \dfrac{\partial}{\partial\tau}
    \begin{bmatrix}
        u(x,s;\tau)\\
        T(\tau)
    \end{bmatrix}
    =\\
    \renewcommand\arraystretch{2.1}
    \begin{bmatrix}
        -\dfrac{1}{T}\dfrac{\partial r}{\partial s}-u\dfrac{\partial r}{\partial x}+\dfrac{\partial^2 r}{\partial x^2}+\dfrac{\partial^4 r}{\partial x^4}\\
        -\displaystyle \frac{1}{T^2} \int_0^1\int_0^L{\frac{\partial u}{\partial s}r}\,\mathrm{d}x\mathrm{d}s
    \end{bmatrix}.
    \renewcommand\arraystretch{1}
\end{multline}

Periodic orbits are dense within the chaotic attractor of the KSE \cite{chaosbook}. Consequently, numerous initial guesses, extracted from a chaotic time series based on close recurrences, converge successfully to a periodic orbit. For a successful convergence, $J\to0$ as $\tau\to\infty$. We consider this to be achieved numerically when the cost function drops below $J = 10^{-12}$. In some cases, however, we find local minima with $J > 0$ as $\tau\to\infty$. Previously, non-zero minima were regarded as failed attempts to find a periodic orbit \cite{azimi2022constructing}. Here, we continue these local minima to determine if they correspond to the ghost of a saddle-node bifurcation.

\begin{figure}
    \centering
    \includegraphics[width=\columnwidth]{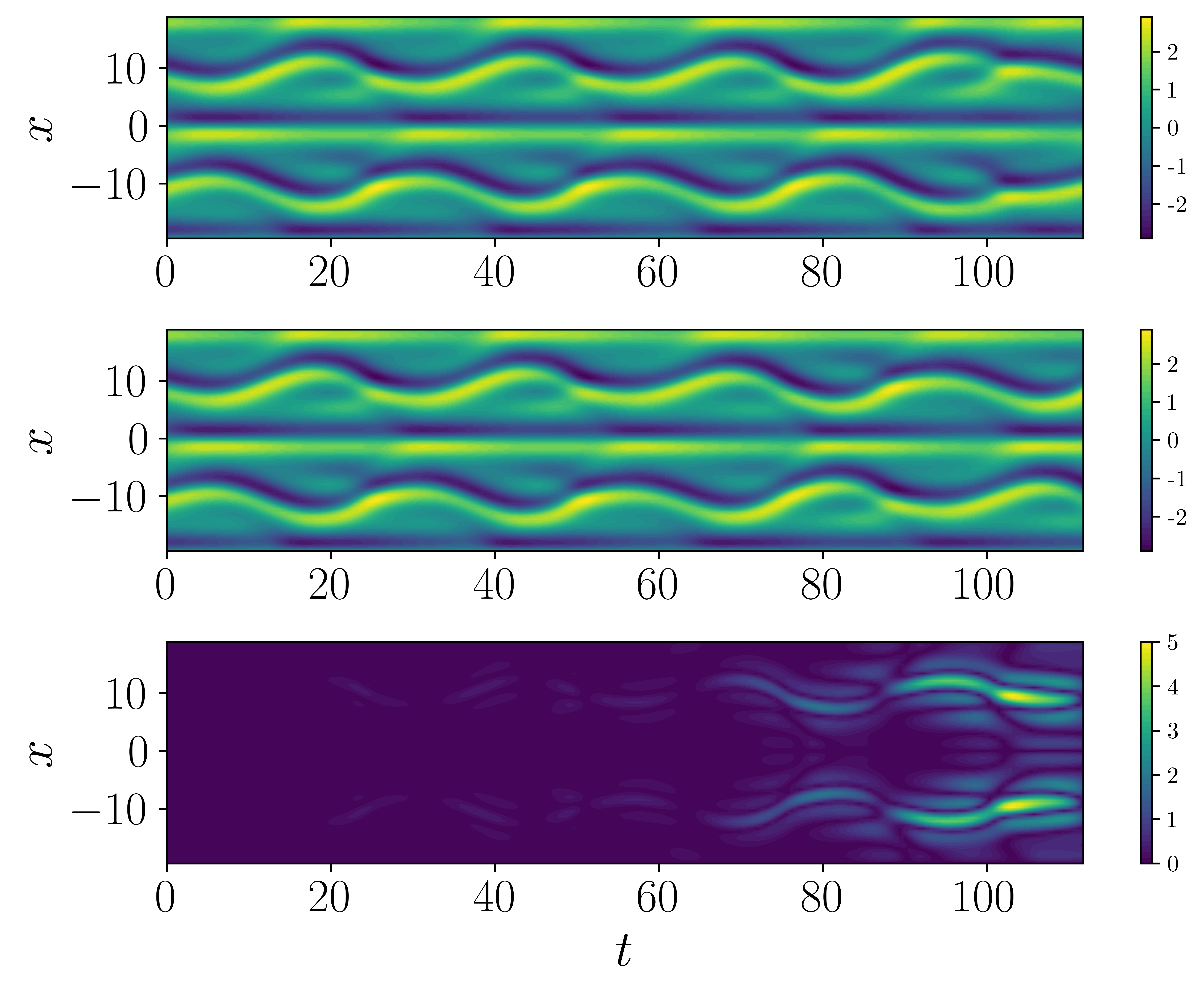}
    \caption{Top: A loop obtained by minimizing the cost function \eqref{eq:KSE_J} for the KSE at $L = 39$ via time-marching the gradient-descent dynamics \eqref{eq:KSE_loop_dynamics}. The convergence saturates around $J \approx 9.72 \times 10^{-5}$ with $T\approx 111.72$. Middle: Time integration of the KSE for a duration $T$, with a point from the loop taken as the initial condition. Bottom: The absolute component-wise difference between the time integration and the loop. The time integration of the KSE follows the loop closely for some time until they diverge as time evolves, suggesting that the loop might be the ghost of a periodic orbit (see Fig.~\ref{fig:KSE_bifurcation}).}
    \label{fig:KSE_triple}
\end{figure}

As an illustration, we study a loop that represents a non-zero local minimum of the cost function, with convergence saturating around $J \approx 9.72 \times 10^{-5}$ with $T\approx 111.72$. Figure \ref{fig:KSE_triple} shows the space--time contours of the loop corresponding to the local minimum of $J$ and the time integration of the dynamics using a point from the loop taken as the initial condition. The trajectory obtained from time integration follows the local minimum for a relatively long time but eventually diverges from it as time evolves. This is evident from the difference between the two, shown in the bottom panel of Fig.~\ref{fig:KSE_triple}. This observation suggests that the loop obtained by minimizing the cost function is probably the ghost of a periodic orbit that does not exist for the chosen parameter value, $L=39$.

\begin{figure}
    \centering
    \includegraphics[width=\columnwidth]{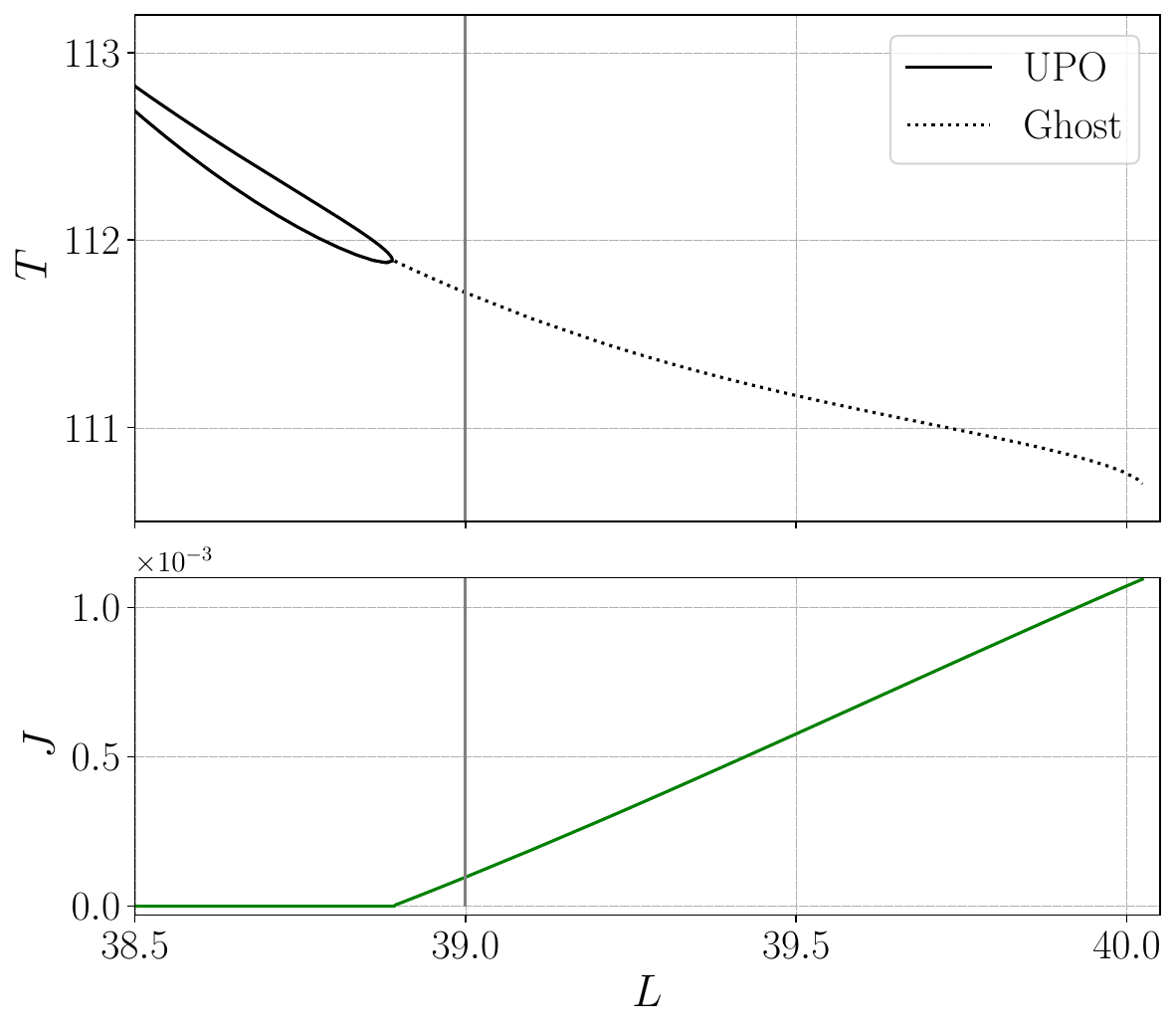}
    \caption{
    Saddle-node bifurcation of an unstable periodic orbit and the resulting ghost for the KSE.
    Top: Bifurcation diagram with the vertical axis showing the period of the structure. The solid and dotted lines represent the unstable periodic solution branch and the ghost branch, respectively.
    Bottom: The cost function \eqref{eq:KSE_J} evaluated for periodic orbits and the ghost.
    The bifurcation occurs at $L^+ \approx 38.889$, below which $J$ reaches the global minimum value of zero. Above the bifurcation value, $J$ takes a local minimum value at the ghost, that lifts from zero as $L$ is increased.
    The vertical line indicates $L=39$, the parameter value for which a local minimum of $J$ was initially computed and then continued (see Fig.~\ref{fig:KSE_triple}).}
    \label{fig:KSE_bifurcation}
\end{figure}

To confirm that the loop corresponding to a local minimum of the cost function is indeed a ghost, we numerically continue its associated solution branch in the direction where the value of the local minimum decreases. By reducing $L$, the cost of the loop eventually decreases to zero at $L^+ \approx 38.889$, indicating that the loop becomes a periodic orbit (see Fig.~\ref{fig:KSE_bifurcation}). At $L=L^+$, two periodic orbits are born in a saddle-node bifurcation. By continuing the two solutions along the parameter $L$, we identify a second saddle-node bifurcation at $L^{-}\approx 36.708$ (not shown), where the two periodic orbits merge and disappear. Therefore, the upper and lower solution branches form an isola in the bifurcation diagram. The ghost branch, namely the contour of $\nabla J = 0$ that is also a minimum, is shown by the dotted line in Fig.~\ref{fig:KSE_bifurcation}. The $\nabla J = 0$ contour also forms a closed curve that passes through the two saddle-node bifurcation points. Above each saddle-node bifurcation, i.e., for $L<L^-$ and $L>L^+$, $\nabla J = 0$ corresponds to a local minimum and thus a ghost, until at some point the type of the extremum changes to a local maximum. 
Inside the isola formed by the periodic solution branch, for $L\in(L^-,L^+)$, the two saddle-node bifurcation points are connected by a contour of $\nabla J = 0$ corresponding to maximum of $J$ (not shown in Fig.~\ref{fig:KSE_bifurcation}).

\begin{figure}
    \centering
    \includegraphics[width=\columnwidth]{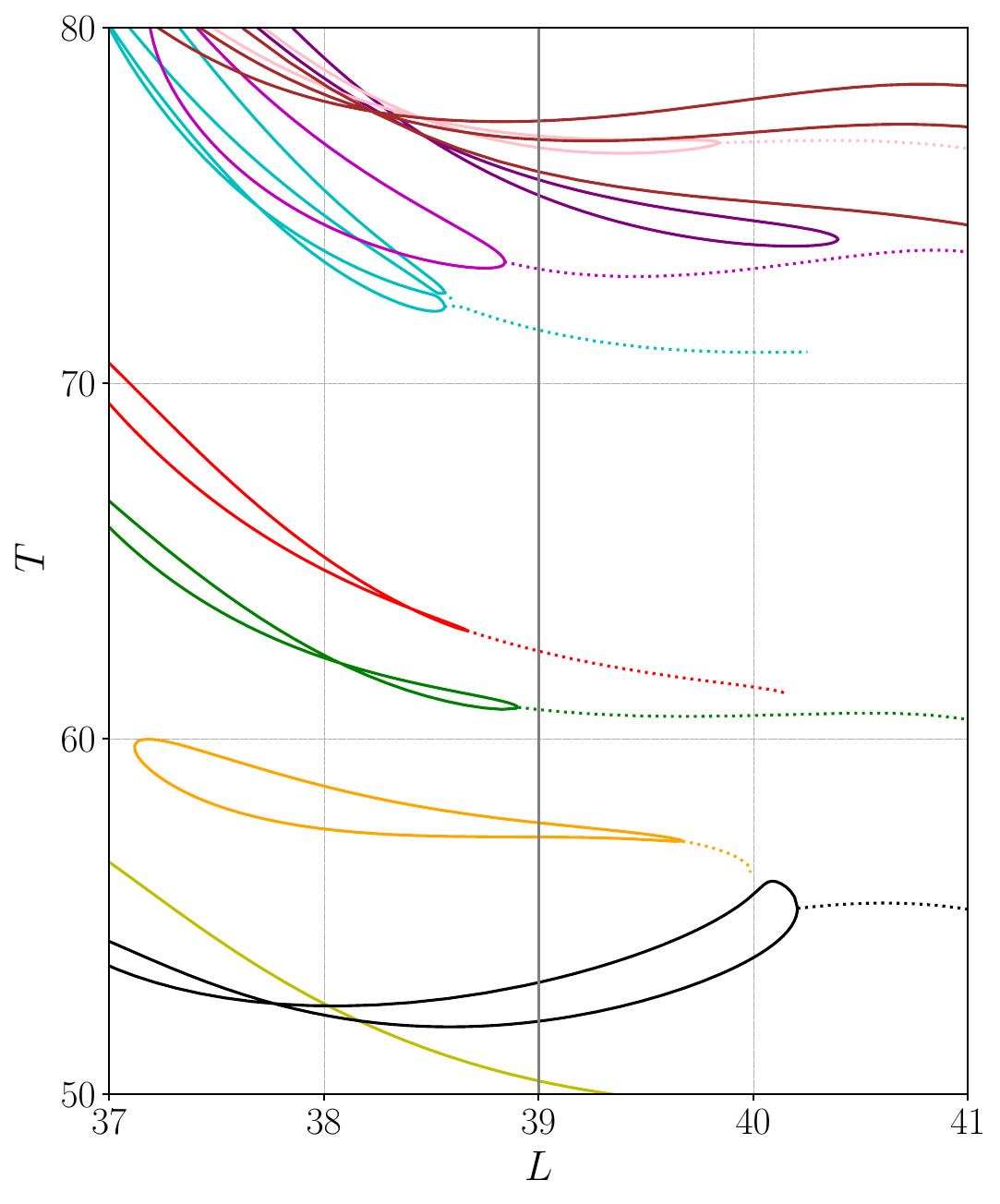}
    \caption{Bifurcation diagram of periodic orbits of the KSE and their ghosts, showing that the state space is `littered' with them. Solid lines represent unstable periodic orbits corresponding to the global minima of the cost function with $\nabla J=0$ and $J=0$. Dotted lines represent ghosts of periodic orbits corresponding to the local minima of the cost function with $\nabla J=0$ and $J>0$. The ghost branches terminate where the extremum of $J$ changes from a local minimum to a local maximum. Note that the ghost of the dark purple bifurcation is very short and hence barely visible. The vertical line indicates $L=39$, the parameter value for which a periodic orbit or a ghost was initially computed and then continued.}
    \label{fig:KSE_litter_plot}
\end{figure}

While we detailed one ghost of a periodic orbit here, we find many such local minima that, upon numerical continuation, are confirmed to be connected to the saddle-node bifurcation of a periodic orbit and thus are ghosts.
Ghosts frequently show up in this chaotic system, and the state space appears to be `littered' with ghosts.
This is exhibited in Fig.~\ref{fig:KSE_litter_plot}, where we pick local minima (or periodic orbits) obtained for $L = 39$ and continue them until we find the related saddle-node bifurcation and periodic orbits (or ghosts).

We observed that the ghost of a periodic orbit is able to shape the trajectories in its vicinity, and that the state space is littered with ghosts of different dynamical relevance. This explains why initial loops, constructed based on close recurrences in a numerical simulation time series, might still fail to converge to a periodic orbit, even when a stringent recurrence criterion is employed \cite{lan2004variational,azimi2022constructing}. Indeed, a close recurrence might emerge as a result of the trajectory following the ghost of a periodic orbit for one cycle. Consequently, the constructed guess is too close to a local minimum of the cost function from which the variational dynamics cannot escape. This motivates the development of new strategies for constructing initial loops \cite{Beck2024} in applications where a large database of periodic orbits at a specific parameter value is needed.
Further, the saddle-node bifurcation connected to the ghost creates two unstable periodic orbits of different periods.
Thus, by continuing ghosts, we are able to find more periodic orbits at other values of the control parameter.

\section{Ghosts in Rayleigh--B\'enard Convection}
\label{sec:ilc}
In this section, we consider the 3D Rayleigh--B\'enard convection (RBC) and demonstrate the relevance of ghosts to the spatial as well as the temporal properties of its dynamics past a saddle-node bifurcation. Specifically, we study the so-called `skewed-varicose' (SV) pattern in the flow within a parameter regime where no equilibrium solution underlying the observed pattern exists. We show that the emergence of this pattern is the result of the dynamics visiting the ghost of an equilibrium solution with a similar pattern.

\subsection{Rayleigh--B\'enard convection and the skewed-varicose pattern}

\begin{figure}
    \centering 
    \includegraphics[width=0.9\columnwidth]{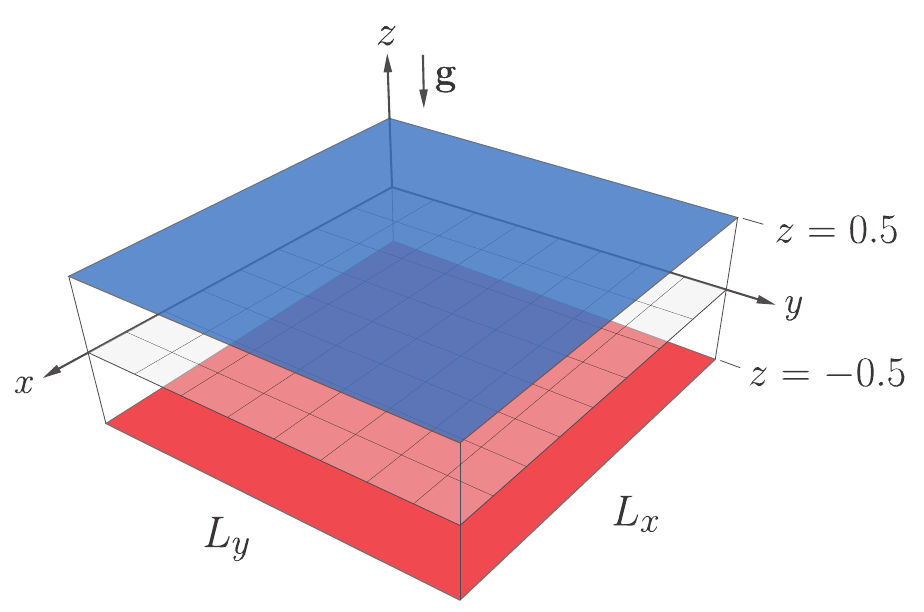}
    \caption{Schematic of the Rayleigh--B\'enard convection cell. The flow is bounded between two parallel walls at $z=\pm0.5$ and is periodic in $x$ and $y$. The walls are kept at constant but different temperatures, such that the confined fluid is heated from the bottom and cooled from the top. The $x$--$y$ plane at $z=0$ is used for flow visualization.} 
    \label{ILC_figure}  
\end{figure}

The 3D RBC describes the flow between two stationary, parallel and horizontal walls that are kept at constant yet different temperatures. The confined fluid is heated from the bottom and cooled from the top, while gravity acts in the wall-normal direction (see Fig.~\ref{ILC_figure}). The Rayleigh ($Ra$) and Prandtl ($Pr$) numbers constitute the two control parameters of the RBC. The former is proportional to the wall temperature difference, thus determining the intensity of thermal driving, while the latter is a fluid property. The RBC is governed by the Oberbeck--Boussinesq equations (OBEs), that in non-dimensional form read
\begin{gather}
    \dfrac{\partial \mathbf{u}_t}{\partial t} = -(\mathbf{u}_t \cdot \nabla)\mathbf{u}_t -\nabla p_t + \sqrt{\frac{Pr}{Ra}} \Delta \mathbf{u}_t + \theta_t\mathbf{e}_z, \label{eq:ghosts_OBEs_momentum}\\
    \dfrac{\partial \theta_t}{\partial t} = -(\mathbf{u}_t \cdot \nabla)\theta_t + \frac{1}{\sqrt{PrRa}} \Delta \theta_t, \label{eq:ghosts_OBEs_energy}\\
    \nabla \cdot \mathbf{u}_t= 0. \label{eq:ghosts_OBEs_continuity}
\end{gather}
Here, $\mathbf{u}_t(\mathbf{x};t) = [u, v, w](x,y,z;t)$, $p_t(\mathbf{x};t)$ and $\theta_t(\mathbf{x};t)$ represent the total velocity vector, pressure and temperature, respectively, at a point $\mathbf{x}$ and time $t$. The unit vector in the positive $z$-direction is indicated with $\mathbf{e}_z$.

We consider RBC in the computational domain depicted in Fig.~\ref{ILC_figure}. In this domain, the OBEs are subject to periodic BCs in the wall-parallel directions $x$ and $y$. At the walls, the flow is subject to no-slip and fixed-temperature BCs:
\begin{gather}
  \mathbf{u}_t(x,y,z=\pm0.5;t)=\mathbf{0},\label{eq:ghosts_RBC_noslip}\\
  \theta_t(x,y,z=\pm0.5;t)=\mp0.5\label{eq:ghosts_RBC_fixed_T}.
\end{gather}
Equations \eqref{eq:ghosts_OBEs_momentum}--\eqref{eq:ghosts_OBEs_continuity} together with the BCs \eqref{eq:ghosts_RBC_noslip} and \eqref{eq:ghosts_RBC_fixed_T} admit the 1D conduction solution: 
\begin{gather}
    \mathbf{u}_0(x,y,z) = \mathbf{0}, \label{eq:ghosts_RBC_laminar_vel} \\
    \theta_0(x,y,z) = -z. \label{eq:ghosts_RBC_laminar_temp}
\end{gather}
We denote the deviation of the total velocity and temperature from the conduction base state with $\mathbf{u} := \mathbf{u}_t - \mathbf{u}_0 =\mathbf{u}_t$ and $\theta := \theta_t - \theta_0$, respectively.

The RBC is equivariant under reflection across the $y=0$ plane,
\begin{equation}
    \pi_y: [u,v,w,\theta](x,y,z) \mapsto [u,-v,w,\theta](x,-y,z),
    \label{eq:ghosts_RBC_sym1}
\end{equation}
reflection across the $y$-axis,
\begin{multline}
    \pi_{xz}:\\ [u,v,w,\theta](x,y,z) \mapsto [-u,v,-w,-\theta](-x, y,-z),
    \label{eq:ghosts_RBC_sym2}
\end{multline}
and continuous translation in the $x$ and $y$ directions,
\begin{multline}
    \tau(\ell_x, \ell_y):\\ [u,v,w,\theta](x,y,z) \mapsto [u,v,w,\theta](x + \ell_x, y + \ell_y,z). \label{eq:ghosts_RBC_sym3}
\end{multline}
These symmetry operations form the equivariance group of the system, which consists of all products of the generators $S_{RBC} \equiv \braket{\pi_y, \pi_{xz}, \tau(\Delta x, \Delta y)}$. The numerical simulations are performed using the open-source software package \textit{Channelflow 2.0} \citep{Gibson2019}. We refer readers to Refs.~\cite{Reetz2020a, reetz2020invariant, Zheng2024part1, Zheng2024part2} for more details on the numerical methods used in this research.

For strong enough thermal driving, the system is driven out of the conduction equilibrium state, giving rise to a large variety of self-organized flow patterns \citep{Cross1993}. Here, we investigate the appearance of one of these patterns, namely the `skewed-varicose' (SV) pattern. This pattern consists of spatially localized distorted convection rolls within otherwise regular straight rolls \citep{Busse1979, Bodenschatz2000}. See panels (b) and (e) of Fig.~\ref{ILC-SV-DNS} for the SV pattern in the temperature field.

In the following, we study the transient emergence of the SV pattern in the dynamical transition from an unstable to a stable equilibrium state. For the chosen parameters, no equilibrium solution capturing the SV pattern exists. We show that the ghost of an equilibrium solution featuring a similar pattern is responsible for the observed flow structure.

\subsection{Equilibria and ghosts underlying the skewed-varicose pattern}
\label{sec:equilibria_ghosts_underlying_SV}

\begin{figure*}
    \centering 
    \includegraphics[width=2.05\columnwidth]{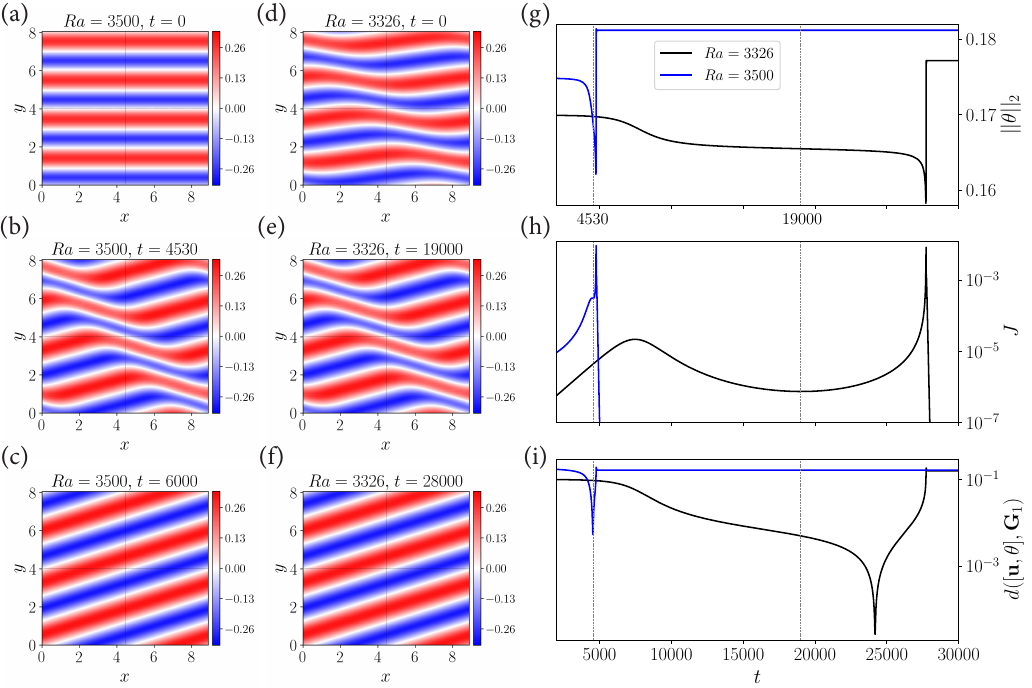}
    \caption{The emergence of the skewed-varicose (SV) pattern during the transition from an unstable to a stable equilibrium solution. This results from the dynamics visiting the ghost of a solution with a similar pattern. Snapshots of the mid-plane temperature in the DNS at $Ra=3500$ are shown in (a)--(c) and at $Ra=3326$ in (d)--(f). The respective initial conditions are marked as open circles in the bifurcation diagram in Fig.~\ref{ILC-SV}. Both flows transiently exhibit the SV pattern ((b) and (e)) and eventually reach the same stable equilibrium ((c) and (f)). Time series of the $L_2$-norm of the temperature deviation from the conduction base state ($\|\theta\|_2$), the state-space velocity ($J$), and the distance from the ghost 1 state $\mathbf{G}_1$ ($d([\mathbf{u},\theta],\mathbf{G}_1)$) are shown in panels (g) to (i), respectively. In (i), the distance from the ghost shows a meaningful drop in both cases, indicating the visiting of the ghost. In (h), a significant slowdown is observed at $Ra=3326$, that is close to the saddle-node bifurcation, and a modest slowdown at $Ra=3500$, that is far from the bifurcation point. The vertical lines mark the instances shown in (b) and (e).} 
    \label{ILC-SV-DNS}  
\end{figure*}

We consider the RBC in a domain with dimensions $[L_x, L_y, L_z] = [8.884, 8.064, 1]$. Following previous studies, the domain size is chosen to be sufficiently large for the SV pattern to form \citep{Subramanian2016, Reetz2020a}. We constrain the dynamics to the subspace of flow fields that are invariant under subgroups of $S_{SV}:=\left<\pi_y\pi_{xz},\tau(L_x/4, -L_y/4)\right>$. Here, $Pr=1.07$ is fixed, and $Ra$ is considered as the control parameter of the system.

At $Ra=3500$, we observe the SV pattern in the transition between two equilibrium states, as shown in Fig.~\ref{ILC-SV-DNS}. The transition starts from an unstable equilibrium solution of straight convection rolls, perturbed along its single unstable eigendirection within the $S_{SV}$-invariant subspace (panel (a)). The unstable solution consists of four pairs of convection rolls with a wavelength $\lambda=2$ that are elongated in the $x$-direction. Following the notation of \citet{Reetz2020a}, we refer to this equilibrium solutions as $R_{\lambda2}$. For $t\to\infty$, the resulting trajectory is attracted to a stable equilibrium solution. The stable solution consists of four pairs of straight rolls with a wavelength $\lambda=2.6$ that are tilted against the $x$-direction (panel (c)). During the transition, the SV pattern emerges transiently in the flow (panel b).

\begin{figure*}
    \centering 
    \includegraphics[width=2.05\columnwidth]{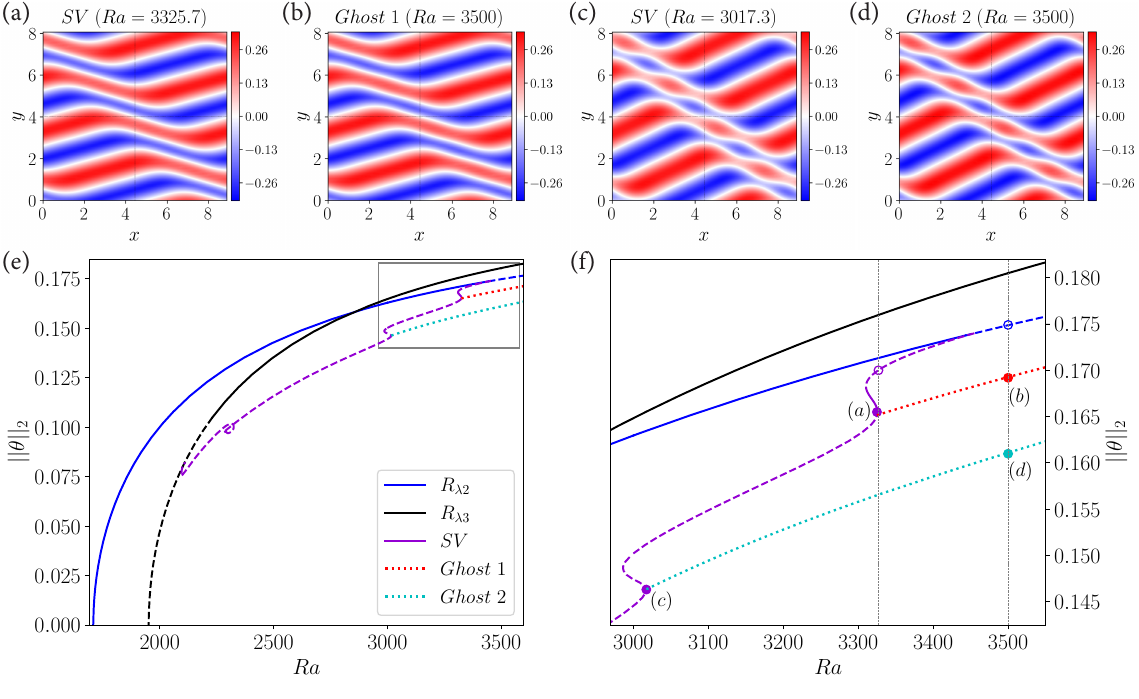}
    \caption{Top row: Mid-plane temperature of the skewed-varicose (SV) equilibrium states at two different saddle-node bifurcations ((a) and (c)) and their ghosts at $Ra=3500$ ((b) and (d)). The appearance of ghost 1 is similar to the SV pattern observed during the transitions at $Ra=3326$ and $3500$, whereas characteristic features of ghost 2 are not observed (see Fig.~\ref{ILC-SV-DNS}). Bottom row: The bifurcation diagram showing the straight-roll equilibria $R_{\lambda_2}$ and $R_{\lambda_3}$, the SV equilibrium solution and the ghost of two of the saddle-node bifurcations. The quantity $\|\theta\|_2$ is the $L_2$-norm of the temperature deviation from the conduction base state. Solid lines indicate stable equilibria, dashed lines denote unstable equilibria, and dotted lines represent the ghost states. Panel (f) shows a zoomed-in view of the region outlined in (e). The location of the snapshots (a)--(d) are marked with filled circles. Open circles indicate the initial conditions used in the DNS performed at $Ra=3326$ (slightly above the saddle-node bifurcation at $Ra=3325.7$) and $Ra=3500$ (see Fig.~\ref{ILC-SV-DNS}). The vertical lines indicate the respective Rayleigh numbers.} 
    \label{ILC-SV}  
\end{figure*}

There is no known equilibrium solution underlying the observed spatial pattern at $Ra = 3500$ \citep{reetz2020invariant}. However, an equilibrium solution branch exhibiting the SV pattern exists at lower Rayleigh numbers, which we refer to as the `SV branch.' The SV branch and other relevant solution branches around the studied parameter value are shown in Fig.~\ref{ILC-SV}e. The SV branch bifurcates at $Ra=2100$ from an equilibrium solution of straight convection rolls with a wavelength $\lambda = 2.8$ that are tilted against the $x$-direction. We refer to this solution branch as $R_{\lambda3}$ following \citet{reetz2020invariant} (see their figure 5a for the temperature field of $R_{\lambda3}$). At $Ra = 3450$, the SV branch joins the $R_{\lambda2}$ equilibrium branch. The SV branch within the $S_{SV}$-invariant subspace is unstable except between the last two saddle-node bifurcations, for $3310.3<Ra<3325.7$. While no SV equilibrium solution exists at $Ra=3500$, the SV pattern emerges transiently in a direct numerical simulation (DNS) at this Rayleigh number. The observed SV pattern at $Ra=3500$ (Fig.~\ref{ILC-SV-DNS}b) appears to be remarkably similar to the SV equilibrium state at the saddle-node bifurcation at $Ra = 3325.7$ (Fig.~\ref{ILC-SV}a). This suggests that the dynamics at $Ra=3500$ visits the ghost of this saddle-node bifurcation.

Unlike the strong similarity in spatial structure, the slowdown associated with the ghost phenomenon is only weakly present in the transition at $Ra=3500$. To illustrate this, we define the state-space velocity for a given velocity field $\mathbf{u}$ and a temperature field $\theta$ as
\begin{equation}
    J^2(\mathbf{u},\theta) := \dfrac{1}{2}\int_\Omega \left[\dfrac{\partial\mathbf{u}}{\partial t}\cdot\dfrac{\partial\mathbf{u}}{\partial t} + \left(\dfrac{\partial\theta}{\partial t}\right)^2\right] \,\mathrm{d}\mathbf{x}.
    \label{eq:ghosts_state_space_velocity}
\end{equation}
Here, $\cdot$ indicates the Euclidean inner product in $\mathbb{R}^3$, $\Omega$ denotes the 3D computational domain, and $\partial_t\mathbf{u}$ and $\partial_t\theta$ are given by the OBEs \eqref{eq:ghosts_OBEs_momentum}--\eqref{eq:ghosts_OBEs_continuity}. We denote the state-space velocity by $J$ as it also serves as the non-negative cost function whose global minima (zeros) correspond to equilibrium solutions, and its minima with a non-zero value can be used to characterize the ghost states, following the methods discussed in Sec.~\ref{sec:method}.

Figure~\ref{ILC-SV-DNS}, panels (g) and (h), shows the DNS time series of the $L_2$-norm of the perturbative temperature $\|\theta\|_2$ and the state-space velocity $J$, respectively. In the DNS at $Ra=3500$, the initial accelerating trend of the state-space velocity changes as the SV pattern emerges. At $t=4530$, when the SV pattern is well captured (panel (b)), the state-space velocity shows a very shallow and brief slowdown before accelerating again. Therefore, although the DNS snapshot is indistinguishably similar to the equilibrium solution at the bifurcation point, the studied $Ra$ is too far from the saddle-node bifurcation for the typical slowdown associated with the ghost phenomenon to occur.

We perform another DNS at $Ra = 3326$, which is only slightly above the saddle-node bifurcation. Consequently, the slowdown due to the ghost is more pronounced in this case. This DNS starts from the unstable SV branch at $Ra=3326$ perturbed along its single unstable eigendirection (Fig.~\ref{ILC-SV-DNS}, panel (d)). The resulting trajectory eventually gets attracted to the same stable equilibrium state as the previous DNS (panel (f)). Note that the straight-roll equilibrium solution $R_{\lambda_2}$ is linearly stable at $Ra = 3326$; hence, a different initial condition from the previous DNS is chosen. Similar to the DNS at $Ra=3500$, the initial accelerating trend of the state-space velocity turns into a decelerating trend when the SV pattern emerges. The state-space velocity reaches a minimum at $t=19000$ when the SV pattern is well captured (panel (e)). The minimum state-space velocity of the DNS at $Ra = 3326$ is more than two orders of magnitude lower than that at $Ra = 3500$. As a result of this slow evolution during the transition, a long plateau is observed in $\|\theta\|_2$ for $Ra = 3326$---the typical temporal characteristic of the ghost phenomenon.

These results suggest that the relevance of a ghost is not limited to situations where the system is extremely close to the saddle-node bifurcation point in the parameter space. The ghost may still play a significant role in the spatial characteristics of the dynamics, even if the control parameter is far from the bifurcation value, and hence the ghost is too weak to induce a significant slowdown in nearby trajectories. In the following, we compute and track the ghost state responsible for the transient appearance of the SV pattern in the flow.

\subsection{Computing and continuing the ghost states}

The RBC has the form of the general dynamical system \eqref{eq:general_governing_equation} with $d=3$ and $n=4$, where pressure is not governed by an explicit evolution equation but adapts itself to the velocity such that the incompressibility constraint is satisfied. The state space of the OBEs contains divergence-free velocity fields $\mathbf{u}$ and temperature fields $\theta$ that, in their perturbative form, satisfy homogeneous Dirichlet BCs at the walls. We employ the state-space velocity $J$, defined in Eq.~\eqref{eq:ghosts_state_space_velocity}, to penalize the deviation of $[\mathbf{u},\theta]$ in the state space of the OBEs from the equilibrium state.

We derive the gradient-descent dynamics for minimizing the cost function within the state space of the OBEs as:
\begin{gather}
	\label{eq:ghosts_RBC_rate_of_u}
    \begin{multlined}
        \dfrac{\partial\mathbf{u}}{\partial\tau} = -\left(\nabla\mathbf{r}_1\right)\mathbf{u} + \left(\nabla\mathbf{u}\right)^\top\mathbf{r}_1-\sqrt{\frac{Pr}{Ra}}\Delta\mathbf{r}_1 \\
        + r_2\nabla\left(\theta + \theta_{0}\right) - \nabla\phi,
    \end{multlined}\\
	\label{eq:ghosts_RBC_rate_of_temp}
    \dfrac{\partial \theta}{\partial \tau} = - \mathbf{r}_1 \cdot \mathbf{e}_z - (\nabla r_2) \cdot \mathbf{u} - \frac{1}{\sqrt{PrRa}}\Delta r_2,\\
    \label{eq:ghosts_RBC_adj_incompressibility}
    \nabla\cdot\mathbf{u} = 0,
\end{gather}
where $\tau$ is the fictitious time parameterizing the gradient-descent dynamics, and $\mathbf{r}_1$ and $r_2$ represent the residuals (i.e., the right-hand side operators) of the momentum equation \eqref{eq:ghosts_OBEs_momentum} and the energy equation \eqref{eq:ghosts_OBEs_energy}, respectively. The gradient-descent dynamics are subject to the following BCs at the walls:
\begin{gather}
    \mathbf{u}(x,y,z=\pm0.5;\tau)=\mathbf{0},\\
    \mathbf{r}_1(x,y,z=\pm0.5;\tau)=\mathbf{0},\\    
    \theta(x,y,z=\pm0.5;\tau)=0,\\
    r_2(x,y,z=\pm0.5;\tau)=0,
\end{gather}
and periodic BCs in the wall-parallel $x$ and $y$ directions. The scalar field $\phi$ is introduced to enforce the incompressibility constraint, analogous to the pressure in the OBEs. The derivation and numerical implementation of the gradient-descent dynamics follow the methodology developed in Ref.~\cite{ashtari2023identifying}. We consider the global minimum (zero) of the cost function, and thus an equilibrium solution, to be achieved numerically if the cost function falls below $J = 10^{-12}$. Otherwise, the converged state is considered to be a non-zero minimum of the cost function, and hence a potential ghost of an equilibrium.

\begin{figure}
    \centering 
    \includegraphics[width=\columnwidth]{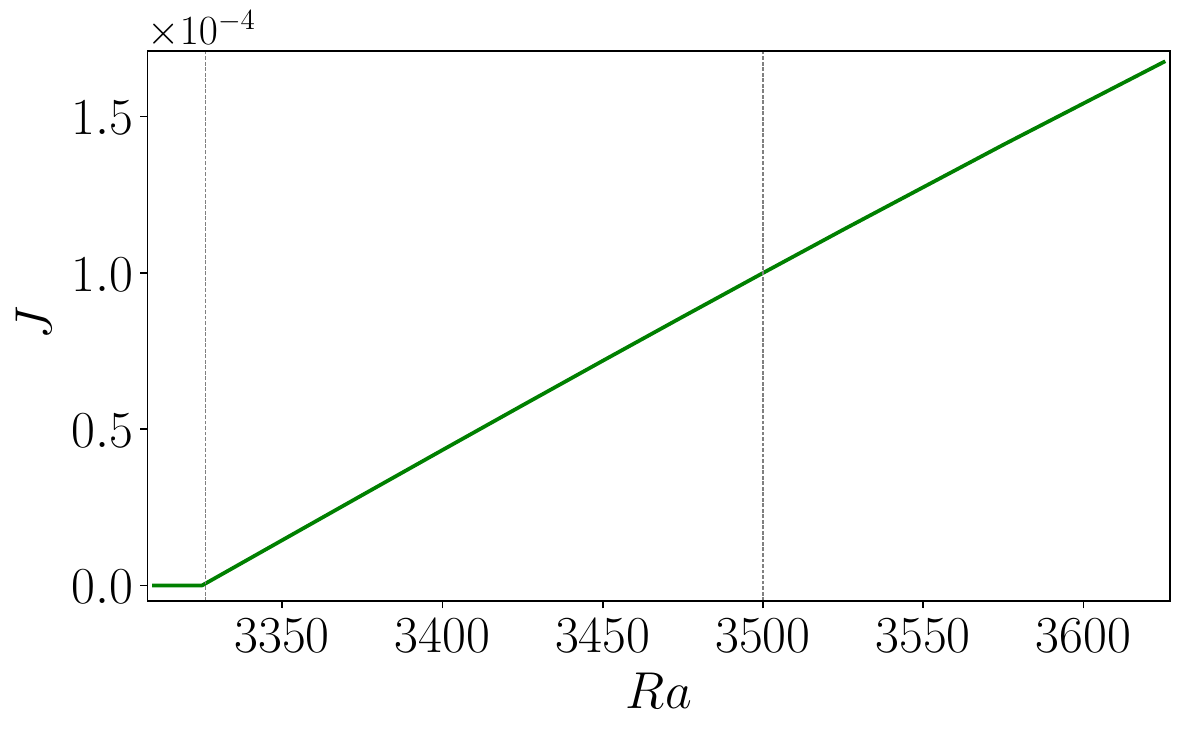}
    \caption{The minimum value of the cost function corresponding to $\tau\to\infty$ in the gradient-descent dynamics \eqref{eq:ghosts_RBC_rate_of_u}--\eqref{eq:ghosts_RBC_adj_incompressibility} as a function of the Rayleigh number. Below the saddle-node bifurcation at $Ra=3325.7$, two equilibrium solutions coexist and thus $J=0$. Past the bifurcation point, the minimum value of the cost function associated with the ghost $1$ lifts from zero. The vertical lines mark $Ra=3326$ and $Ra = 3500$ at which DNS were performed (see Fig.~\ref{ILC-SV-DNS}).} 
    \label{ILC-BD-J}  
\end{figure}

We compute the ghost branch originating from the saddle-node bifurcation of the SV branch at $Ra=3325.7$ and label it as the \textit{ghost 1} branch. This branch is presented in the bifurcation diagram shown in Fig.~\ref{ILC-SV}. The cost function $J$ is zero for Rayleigh numbers at and below the bifurcation value and increases with the Rayleigh number above it, as shown in Fig.~\ref{ILC-BD-J}. To characterize how the dynamics visit the ghost, we quantify the distance of the evolving flow fields from the ghost state. We define the distance in the state space of the OBEs as:
\begin{multline}
    d^2\left(
    \begin{bmatrix}    
        \mathbf{u}_1\\
        \theta_1
    \end{bmatrix},\begin{bmatrix}
        \mathbf{u}_2\\
        \theta_2
    \end{bmatrix}\right) := \\
    \int_\Omega\left[\left(\mathbf{u}_1-\mathbf{u}_2\right)\cdot\left(\mathbf{u}_1-\mathbf{u}_2\right) + (\theta_1-\theta_2)^2\right]\mathrm{d}\mathbf{x}.
    \label{eq:ghosts_distance_RBC}
\end{multline}
The DNS time series of the distance between the instantaneous flow fields and the respective ghost 1 at $Ra=3500$ and $Ra=3326$ is shown in Fig.~\ref{ILC-SV-DNS}i. For both Rayleigh numbers, the visit of the ghost state by the trajectory is reflected in a significant drop in the distance from the ghost. For $Ra=3500$, the minimum distance from the respective ghost state coincides closely with the minimum state-space velocity during the slowdown episode. For $Ra=3326$, however, the trajectory reaches its minimum distance from the ghost state with a delay compared to the moment when the minimum state-space velocity is achieved. This happens because the level sets of the state-space velocity $J$ \eqref{eq:ghosts_state_space_velocity} are not hyperspheres around the ghost state, given the distance \eqref{eq:ghosts_distance_RBC} as the metric. Therefore, when the minimum state-space velocity is observed as a result of the trajectory becoming tangent to a level set of $J$, the trajectory is not tangent to a hypersphere of constant distance from the ghost state. Consequently, the minimum distance from the ghost and the minimum state-space velocity are not achieved at the same time. They take their minimum values simultaneously only when the trajectory passes through the ghost state itself.

The SV branch undergoes multiple other saddle-node bifurcations, each leaving its own ghost in the state space. We compute the ghost branch originating from the saddle-node bifurcation at $Ra=3017.3$ and label it as the \textit{ghost 2} branch (see the bifurcation diagram in Fig.~\ref{ILC-SV}). The ghost states on this branch have a different pattern structure along the domain diagonal compared to ghost 1 (see panels (c) and (d) of Fig.~\ref{ILC-SV}, compare with panels (a) and (b)). Similar to ghost 1, we observe that by increasing $Ra$ from the bifurcation value, the spatial pattern of the ghost states remains very similar to that of the parent solution at the saddle-node bifurcation. However, the ghost loses its potential to induce a significant slowdown as the associated state-space velocity increases.

The minimum distance of the DNS trajectories from ghost 2 (not shown in Fig.~\ref{ILC-SV-DNS}) is three orders of magnitude larger than the minimum distance from ghost 1 at $Ra=3326$, and one order of magnitude larger at $Ra=3500$. As a result, the ghost 2 states, and similarly the ghosts of the rest of the saddle-node bifurcations, do not influence the temporal behavior of the dynamics compared to ghost 1 in the particular simulations performed in Sec. \ref{sec:equilibria_ghosts_underlying_SV}. Unsurprisingly, the characteristic diagonal features of ghost 2 do not appear in the instantaneous flow fields in the DNS.

\section{Summary and concluding remarks}
\label{sec:conclusion}
The phase portrait of a dynamical system is organized by dynamically connected time-invariant sets such as unstable equilibria and periodic orbits. Whenever two such solutions collide and disappear through a saddle-node bifurcation, their properties continue to be felt in the region of the state space where the solutions used to exist, a phenomenon known as the ghost of the bifurcation. These properties include slow evolution in the case of bifurcating equilibria and near-periodic behavior in the case of bifurcating periodic orbits. Close to the saddle-node bifurcation point in the parameter space, these properties are felt strongly in the ghost region of the state space, but they lose their significance as the control parameter of the system is varied farther from the bifurcation point.

In this paper, we have presented different examples of quantitatively characterizing ghosts near saddle-node bifurcations in discrete- and continuous-time dynamical systems. These examples cover both fixed points and periodic orbits, but there is no reason the presented methods could not be extended to the saddle-node bifurcations of more exotic invariant solutions, such as invariant tori \citep{parker2023}. Our analyses are based on formalizing the ghost phenomenon by defining representative state-space structures for this phenomenon, which we refer to as the ghost states. For all types of invariant solutions, the key is to define a non-negative cost function for state-space sets with prescribed topology, which is zero if and only if evaluated on an invariant solution. This formulation allows us to define the ghost states as non-zero minima of the cost function that appear as a result of the saddle-node bifurcation. Our definition does not require the system to be asymptotically close to the saddle-node bifurcation point in the parameter space. Instead, it enables us to track the ghost properties as the system moves away to finite distances from the bifurcation point. If such methods are computationally feasible for finding invariant solutions, as has been demonstrated \citep{parker2022variational, ashtari2023identifying}, then they will also be feasible for finding the ghost states.

Characterizing the ghost phenomenon in terms of the ghost states and computing them using the methods presented in Sec. \ref{sec:method} has enabled two principal contributions. First, it has enabled identifying the spatial characteristics in addition to the temporal characteristics of the ghost phenomenon in spatially extended systems, as the proposed numerical method scales appropriately to high-dimensional discretizations of PDEs. We demonstrated this for the 3D Rayleigh--B\'enard convection, where the dynamics visit the ghost of an equilibrium solution and, consequently, slows down and exhibits a specific spatial pattern structure.

The second contribution is a global characterization of the ghosts that result from the saddle-node bifurcation of time-varying invariant solutions. In this paper, we have considered the specific case of periodic orbits undergoing a saddle-node bifurcation. We computed the ghost of a stable periodic orbit in the Lorenz dynamics, after whose bifurcation the system exhibits intermittent chaos, continuing to return to the vicinity of the ghost. We also demonstrated the method by computing several ghosts of periodic orbits for the 1D Kuramoto--Sivashinsky PDE in a chaotic regime.

We have shown that ghosts are not uniquely defined structures, as their precise details depend on the choice of the cost function. However, close to the bifurcation point in the parameter space, where ghosts are most relevant, such differences were small for the examples studied. Further from the bifurcation point, the differences become more significant, but ghosts are also less relevant for the dynamics.

All the examples presented in this paper are governed by deterministic evolution equations. However, the ghost phenomenon is robust to intrinsic noise and is present in stochastically forced systems as well \citep{sardanyes2020noise}. We leave the quantitative characterization of the ghost phenomenon in terms of the ghost states in such systems for future research.

\begin{acknowledgments}
The authors are thankful to Florian Reetz for insightful discussions. This work was supported by the European Research Council (ERC) under the European Union's Horizon 2020 research and innovation programme (grant no. 865677).
\end{acknowledgments}

% \bibliography{references}


\end{document}